\documentclass[11pt,reqno]{amsart}
\usepackage{amsmath,amsfonts, amssymb, amsthm}
  \usepackage{paralist}
  \usepackage{graphics} 
  \usepackage{epsfig} 
\usepackage{graphicx}  \usepackage{epstopdf}
 \usepackage[colorlinks=true]{hyperref}
\hypersetup{urlcolor=blue, citecolor=red}

\newtheorem{theorem}{Theorem}[section]

\newtheorem{lemma}[theorem]{Lemma}

\newtheorem{conjecture}{Conjecture}[section]

\theoremstyle{definition}

\newtheorem{remark}{Remark}[section]

\def\pmod #1{\ ({\rm{mod}}\ #1)}
\def\Z{\Bbb Z}
\def\N{\Bbb N}

\def\Q{\Bbb Q}

\def\R{\Bbb R}

\def\l{\left}
\def\r{\right}
\def\bg{\bigg}
\def\({\bg(}
\def\){\bg)}
\def\t{\text}
\def\f{\frac}
\def\mo{{\rm{mod}\ }}
\def\pmod#1{\ (\mo\ #1)}

\def\ls{\leq}
\def\gs{\geq}

\def\sm{\setminus}

\def\al{\alpha}

\def\eq{\equiv}

\def\da{\delta}

\def\Proof{\noindent{\it Proof}}

\begin{document}
\hbox{Preprint, {\tt arXiv:2010.05775}}
\medskip

\title[Sums of four rational squares with certain restrictions]
      {Sums of four rational squares \\with certain restrictions}
\author[Zhi-Wei Sun]{Zhi-Wei Sun}


\address{Department of Mathematics, Nanjing
University, Nanjing 210093, People's Republic of China}
\email{zwsun@nju.edu.cn}

\keywords{Sums of rational squares, Lagrange's theorem.
\newline \indent Supported by the Natural Science Foundation of China (grant no. 11971222).}

\begin{abstract}
In this paper we mainly study sums of four rational squares with certain restrictions.
Let $\mathbb Q_{\ge0}$ be the set of nonnegative rational numbers.
We establish the following four-square theorem for rational numbers:
For any $a,b,c,d\in\mathbb Q_{\ge0}$,
each $r\in\mathbb Q_{\ge0}$ can be written as $x^2+y^2+z^2+w^2$
with $x,y,z,w\in\mathbb Q_{\ge0}$ such that $ax+by+cz+dw$ is a rational square
(or a rational cube). This paper also contains many conjectures; for example,  for any positive integers $a$ and $b$ with $\gcd(a,b)=1$, we conjecture that each $r\in\Q_{\gs0}$ can be written as
$aw^4+bx^4+y^2+z^2$ with $w,x,y,z\in\Q$.
\end{abstract}
\maketitle

\section{Introduction}

Lagrange's four-square theorem established in 1770 states that each $n\in\N=\{0,1,2,\ldots\}$ can be written as $x^2+y^2+z^2+w^2$ with $x,y,z,w\in\N$. In 2017 the author \cite{JNT135} proved that
this can be refined in various ways.

As in \cite{JNT135}, we call a polynomial
$P(x,y,z,w)\in\Z[x,y,z,w]$ {\it suitable} if any $n\in\N$ can be written as
$x^2+y^2+z^2+w^2$ with $x,y,z,w\in\N$ such that $P(x,y,z,w)$ is a square.
The author \cite{JNT135} showed that the linear polynomial $x,2x,x-y,2(x-y)$ are suitable, and conjectured that
$$x+2y,\, x+3y,\ x+24y,\ 2x-y,\ 4x-3y,\ 6x-2y$$
are also suitable. Based on the idea of \cite{JNT135}, Y.-C. Sun and the author \cite{SunSun} proved that $x+2y$ is suitable, moreover any $m\in\Z^+=\{1,2,3,\ldots\}$ can be written as $x^2+y^2+z^2+w^2$
with $x,y,z,w\in\N$ such that $x+2y$ is a positive square. Recently, Y.-F. Sue and H.-L. Wu \cite{SW20}
proved that $x+3y$ is also suitable via the arithmetic theory of ternary quadratic forms, and A. Machiavelo et al. \cite{MT,MRT} proved
the author's 1-3-5 conjecture which states that
$x+3y+5z$ is suitable by using Hamilton quaternions.

For convenience, we let $\Q_{\ge0}$ denote the set of nonnegative rational numbers.
For $h=2,3$ we set
$$\Z^h=\{t^h:\ t\in\Z\}\ \ \t{and}\ \ \Q^h=\{t^h:\ t\in\Q\}.$$
For $a,b,c,d\in\Z$ and $r=m/n\in\Q_{\gs0}$ with $m\in\N$ and $n\in\Z^+=\{1,2,3,\ldots\}$, if $mn^3=x^2+y^2+z^2+w^2$
with $x,y,z,w\in\N$ and $ax+by+cz+dw\in\Z^2$, then
$$r=\f{mn^3}{n^4}=\l(\f x{n^2}\r)^2+\l(\f y{n^2}\r)^2+\l(\f{z}{n^2}\r)^2+\l(\f{w}{n^2}\r)^2$$
and
$$a\f x{n^2}+b\f y{n^2}+c\f z{n^2}+d\f w{n^2}\in\Q^2.$$
Thus the confirmed 1-3-5 conjecture implies that each $r\in\Q_{\gs0}$
can be written as $x^2+y^2+z^2+w^2\ (x,y,z,w\in\Q_{\gs0})$ such that $x+3y+5z$
is a rational square.

In this paper we study sums of
the form $x^2+y^2+z^2+w^2\ (x,y,z,w\in\Q_{\gs0})$
with $ax+by+cz+dw$ a rational square or a rational cube, where
$a,b,c,d\in\Q_{\gs0}$.

Now we state our results in this paper.

\begin{theorem}  [Four-square Theorem for Rational Numbers]
\label{Th-abcd} Let $a,b,c,d$ and $h\in\{2,3\}$. For each $r\in\Q_{\gs0}$, there are $x,y,z,w\in\Q_{\gs0}$ such that
$r=x^2+y^2+z^2+w^2$ and $ax+by+cz+dw\in \Q^h$.
Moreover, there is a positive integer $N$ depending only on $a,b,c,d$
such that any $r\in\N$ can be written as $x^2+y^2+z^2+w^2$ with $x,y,z,w\in\N/N$
and $ax+by+cz+dw\in\Q^h$.
\end{theorem}

For $h\in\{2,3\}$ and $a,b,c,d\in\Q_{\gs0}$, it is interesting to find the exact value
of $S_h(a,b,c,d)$, the smallest positive integer $N$ such that any $m\in\N$ can be written as $x^2+y^2+z^2+w^2$ with $x,y,z,w\in\N/N$
and $ax+by+cz+dw\in\Q^h$.

\begin{theorem}  \label{x+y-square} Any $m\in\N$ can be written as $x^2+y^2+z^2+w^2$
with $x,y,z\in\N/3$ and $w\in\N$ such that $x+y$ is an integer square.
\end{theorem}
\begin{remark} The author \cite{JNT135} proved that any $m\in\N$ can be written as
$x^2+y^2+z^2+w^2$
with $x,y,z,w\in\N$ such that $x-y$ is an integer square.
\end{remark}

\begin{theorem}  \label{x+y+z-square} Any $m\in\N$ can be written as $x^2+y^2+z^2+w^2$
with $x,y,z,w\in\N/27$ such that $x+y+z$ is an integer square.
\end{theorem}
\begin{remark} The author \cite{JNT135} proved that any $m\in\N$ can be written as
$x^2+y^2+z^2+w^2$
with $x,y,z,w\in\Z$ such that $x+y+z$ is an integer square.
\end{remark}

\begin{theorem} \label{four-square}
Any $m\in\N$ can be written as $x^2+y^2+z^2+w^2$
with $x,y,z,w\in\N/50$ such that $x+y+z+w$ is a rational square.
\end{theorem}
\begin{remark} Y.-C. Sun and the author \cite{SunSun} proved that any $m\in\N$ can be written as
$x^2+y^2+z^2+w^2$
with $x,y,z,w\in\Z$ such that $x+y+z+w$ is an integer square.
\end{remark}

\begin{theorem}\label{Th-112} Any $m\in\N$ can be written as $x^2+y^2+z^2+w^2$
with $x\in\N$ and $y,z,w\in\N/3$ such that $x+y+2z$ is an integer square.
\end{theorem}
\begin{remark}\label{Rem1.1} By \cite[Theorem 1.4]{S19IJNT}, any $m\in\Z^+$ can be written
$x^2+y^2+z^2+w^2$ with $x,y,z,w\in\Z$ and $x+y+2z\in\{4^a:\ a\in\N\}$.
\end{remark}

\begin{theorem}\label{Th122} Any $m\in\N$ can be written as $x^2+y^2+z^2+w^2$
with $x,y,z,w\in\N/13$ such that $x+2y+2z$ is an integer square.
\end{theorem}
\begin{remark}\label{Rem1.1} By \cite[Theorem 1.4]{S19IJNT}, any $m\in\Z^+$ can be written
$x^2+y^2+z^2+w^2$ with $x,y,z,w\in\Z$ and $x+2y+2z\in\{4^a:\ a\in\N\}$.
\end{remark}

\begin{theorem}\label{Th123} {\rm (i)} Any $m\in\N$ can be written as $x^2+y^2+z^2+w^2$
with $x,y,z\in\N/3$ and $w\in\N$ such that $x+y+3z$ is an integer square.

{\rm (ii)} Any $m\in\N$ can be written as $x^2+y^2+z^2+w^2$
with $x,y,z,w\in\N/5$ such that $x+2y+3z$ is an integer square.
\end{theorem}
\begin{remark}\label{Rem1.1}  By
\cite[Theorem 1.6(ii) and Theorem 1.7(iv)]{SunSun}, for each $d=1,2$,
any $n\in\N$ can be written as $x^2+y^2+z^2+w^2\ (x,y,z,w\in\Z)$ with $x+dy+3z$ a square.
\end{remark}

\begin{theorem}\label{Th1.2} Any integer $m>1$ with $m\not=10$ can be written as $x^2+y^2+z^2+w^2$ with $x,y,z,w\in\Z$
and $x+2y+3z\in\{4^a:\ a\in\Z^+\}$.
\end{theorem}
\begin{remark}\label{R-Th1.2} In \cite[Conjecture 4.5(ii)]{S19IJNT}, the author conjectured that
any $m\in\Z^+$ can be written as $x^2+y^2+z^2+w^2$ with $x,y,z,w\in\N$
and $|x+2y-3z|\in\{4^a:\ a\in\N\}$.
\end{remark}

\begin{theorem} \label{Th-1123} {\rm (i)} Any $m\in\N$ can be written as $x^2+y^2+z^2+w^2$
with $x,y,z\in\Q_{\gs0}$ and $w\in\N$ such that $x+y+z+2w\in\Q^2$.

{\rm (ii)} Any $m\in\N$ can be written as $x^2+y^2+z^2+w^2$
with $x,y,z,w\in\N/5$ such that $x+y+2z+3w\in\Z^2$.

{\rm (iii)} Any $m\in\N$ can be written as $x^2+y^2+z^2+w^2$
with $y\in\N$ and $x,z,w\in\Q_{\gs0}$ such that $x+2y+3z+3w\in\Z^2$.
\end{theorem}

We are going to prove Theorem 1.1 in the next section.
Section 3 is devoted to our proofs of Theorems 1.2-1.4.
Theorems \ref{Th-112}-\ref{Th122}
and Theorems \ref{Th123}-\ref{Th-1123} will be proved in Sections 4 and 5 respectively.
We will pose some conjectures in Section 6.

\section{Proof of Theorem \ref{Th-abcd}}

\setcounter{equation}{0}
 \setcounter{conjecture}{0}
 \setcounter{theorem}{0}

 The Gauss-Legendre theorem (cf. \cite[p.\,23]{N96}) states that $n\in\N$ can be written as the sum of three squares if and only if $n$ does not belong to the set
 \begin{equation}\label{E}E=\{4^s(8t+7):\ s,t\in\N\}.
 \end{equation}

 \begin{lemma} \label{mn} Let $m,n\in\N$.

 {\rm (i)} Assume $16\nmid m$, and  suppose that
 $n$ is odd if and only if $m\eq2\pmod4$ or $m\eq 4,7\pmod 8$. Then $m-n^4,m-n^6\not\in E$.

 {\rm (ii)} Assume $m=16m_0$ with $m_0$ odd, and suppose that $4\mid n$ if $m_0\not\eq7\pmod8$,
 and $n\eq2\pmod4$ otherwise. Then $m-n^4,m-n^6\not\in E$.

 {\rm (iii)} Suppose that $m\eq32\pmod{64}$. Then $m-n^4\not\in E$ if $4\mid n$,
  and $m-n^6\not\in E$ if $2\mid n$.

 \end{lemma}
 \Proof.  If $m\eq2\pmod 4$ or $m\eq4,7\pmod 8$, and $n$ is odd,
 then
 $$m-n^4\eq m-n^6\eq m-1\eq 1,3,5,6\pmod8$$
 and hence $m-n^4,m-n^6\not\in E$.
 If $2\nmid m$, $m\not\eq7\pmod8$, and $n$ is even, then
 $$m-n^4\eq m-n^6\eq m\eq1,3,5\pmod 8$$
 and hence $m-n^4\not\in E$.
 If $m\eq8\pmod{16}$, and $n$ is even, then
 $$m-n^4\eq m-n^6\eq m\eq 8\pmod{16}$$
 and hence $m-n^4,m-n^6\not\in E$.

Now we consider the case $m=16m_0$ with $m_0$ odd.
If $m_0\not\eq7\pmod8$ and $n\eq0\pmod4$, then
$$m_0-\l(\f n2\r)^4\eq m_0-4\l(\f n2\r)^6\eq m_0\not\eq7\pmod 8$$
and hence $m-n^4,m-n^6\not\in E$.
If $m_0\eq7\pmod8$ and $n\eq2\pmod4$, then
$$m_0-\l(\f n2\r)^4\eq7-1=6\pmod 8$$
and
$$m_0-4\l(\f n2\r)^6\eq 7-4\eq3\pmod 8,$$
therefore $m-n^4,m-n^6\not\in E$.

Now assume that $m\eq32\pmod{64}$.  If $4\mid n$, then $m-n^4\eq32\pmod 64$ and hence $m-n^4\not\in E$.
If $2\mid n$, then $m-n^6\eq 32\pmod {64}$ and hence $m-n^6\not\in E$.

In view of the above discussion, we have proved the desired result. \qed

\begin{lemma}\label{Lem-SS} Let $a,b,c,d\in\N$ with $N=a^2+b^2+c^2+d^2>0$,
and let $h\in\{2,3\}$ and $m,n\in\N$ with $mN-n^{2h}\in \N\sm E$.
Then there are $x,y,z,w\in\Z/N$
satisfying
 \begin{equation}\label{mn^h}\begin{cases}x^2+y^2+z^2+w^2=m,&
 \\ax+by+cz+dw=n^h.&
 \end{cases}
 \end{equation}
\end{lemma}
\Proof. We utilize the basic idea in the proof of \cite[Theorem 1.4]{SunSun}.
Let $s=n^h$. As $mN-s^2\in\N\sm E$, by the Gauss-Legendre theorem on sums of three squares,
there are $t,u,v\in\Z$ such that $mN-s^2=t^2+u^2+v^2$. Define
\begin{equation*}\begin{cases} x=\f{as-bt-cu-dv}{N},
\\y=\f{bs+at+du-cv}{N},
\\z=\f{cs-dt+au+bv}{N},
\\w=\f{ds+ct-bu+av}{N}.
\end{cases}\end{equation*}
Then $x,y,z,w\in\Z/N$, and also
\begin{equation*}\begin{cases} ax+by+cz+dw=s,
\\ay-bx+cw-dz=t,
\\az-bw-cx+dy=u,
\\aw+bz-cy-dx=v.
\end{cases}
\end{equation*}
Applying Euler's four-square identity, we obtain
$$N(x^2+y^2+z^2+w^2)=s^2+t^2+u^2+v^2=mN$$
and hence $x^2+y^2+z^2+w^2=m$. This concludes the proof. \qed

 For any real numbers $a_1,\ldots,a_n,x_1,\ldots,x_n$, we have the
Cauchy-Schwarz inequality (cf. \cite[p.\,178]{N96})
$$(a_1x_1+\ldots+a_nx_n)^2\ls (a_1^2+\ldots+a_n^2)(x_1^2+\ldots x_n^2).$$
We will make use of the inequality in our proof of the following lemma.

 \begin{lemma}\label{Lem2.1} Let $a,b,c,d,m$ be nonnegative real numbers with $a^2+b^2+c^2+d^2\not=0$. Suppose that
 $x,y,z,w$ are real numbers satisfying
 \begin{equation}\begin{cases}x^2+y^2+z^2+w^2=m,&
 \\ax+by+cz+dw=s,&
 \end{cases}
 \end{equation}
 where
 \begin{equation}\label{s-ineq}s\gs\sqrt{m(a^2+b^2+c^2+d^2-\min(\{a^2,b^2,c^2,d^2\}\sm\{0\}))}.
 \end{equation}
 Then all the numbers $ax,by,cz,dw$ are nonnegative.
 \end{lemma}
 \Proof. Let
 $$t=ay-bx+cw-dz,\ u=az-bw-cx+dy,\ v=aw+bz-cy-dx.$$
 By Euler's four-square identity, we have
\begin{equation}\label{Euler}(a^2+b^2+c^2+d^2)(x^2+y^2+z^2+w^2)=s^2+t^2+u^2+v^2.\end{equation}
 Since
 \begin{equation}\label{stuv}\begin{cases} ax+by+cz+dw=s,
\\ay-bx+cw-dz=t,
\\az-bw-cx+dy=u,
\\aw+bz-cy-dx=v,
\end{cases}
\end{equation}
we get
\begin{equation}\label{abcd}\begin{cases} x=\f{as-bt-cu-dv}{a^2+b^2+c^2+d^2},
\\y=\f{bs+at+du-cv}{a^2+b^2+c^2+d^2},
\\z=\f{cs-dt+au+bv}{a^2+b^2+c^2+d^2},
\\w=\f{ds+ct-bu+av}{a^2+b^2+c^2+d^2}.
\end{cases}\end{equation}

Suppose that $a>0$. Then
$$s^2\gs m(a^2+b^2+c^2+d^2-a^2)=(b^2+c^2+d^2)m$$
and hence
$$(a^2+b^2+c^2+d^2)s^2\gs(b^2+c^2+d^2)(a^2+b^2+c^2+d^2)m=(b^2+c^2+d^2)(s^2+t^2+u^2+v^2).$$
Thus $a^2s^2\gs (b^2+c^2+d^2)(t^2+u^2+v^2)$.
By the Cauchy-Schwarz inequality,
$$(bt+cu+dv)^2\ls(b^2+c^2+d^2)(t^2+u^2+v^2).$$
Therefore $as\gs|bt+cu+dv|$ and hence $x>0$ in view of \eqref{abcd}.

Similarly, $y\gs0$ if $b>0$, and $z\gs 0$ if $d>0$. This concludes the proof. \qed
\medskip

Now we present an auxiliary theorem.

 \begin{theorem} \label{Th-aux} Let $a,b,c,d\in\N$ with $N=a^2+b^2+c^2+d^2>0$.
 Let $h\in\{2,3\}$ and
 \begin{equation}\label{M_h}M_h(a,b,c,d):=\l(\f{2^{\lfloor\nu_2(N)/(2h)\rfloor+h}}{\root{2h}\of{N}
 -\root{2h}\of{N-\min(\{a^2,b^2,c^2,d^2\}\sm\{0\})}}\r)^{2h},
 \end{equation}
 where $\nu_2(N)=\max\{k\in\N:\ 2^k\mid N\}$ is the $2$-adic valuation of $N$.
 If $m\gs M_h(a,b,c,d)$ is an integer not divisible by $4^h$,
 then there are $x,y,z,w\in\N/N$ such that $x^2+y^2+z^2+w^2=m$
 and $ax+by+cz+dw\in\Z^h$.
 \end{theorem}
 \Proof. Write $N=2^{2hq+r}N_0$ with $q\in\N$, $r\in\{0,\ldots,2h-1\}$ and $N_0\in\{1,3,5,\ldots\}$.
 Then $q=\lfloor\nu_2(N)/(2h)\rfloor$. Suppose that $m=2^sm_0\gs M_h(a,b,c,d)$,
 where $s\in\{0,\ldots,2h-1\}$ and $m_0\in\{1,3,5,\ldots\}$. Clearly
 $mN=2^{2hq+r+s}m_0N_0$.  Let
 $$q'=q+\l\lfloor\f{r+s}{2h}\r\rfloor=\begin{cases}q+1&\t{if}\ r+s\gs 2h,
 \\q&\t{if}\ r+s<2h,\end{cases}$$
 and let $t$ be the least nonnegative residue of $r+s$ modulo $2h$.
 Then $mN=4^{hq'}(2^tm_0N_0)$ with $4^h\nmid 2^tm_0N_0$.
 As $m\gs M_h(a,b,c,d)$, the length of the interval
 $$I=\l[\root{2h}\of{m(N-\min(\{a^2,b^2,c^2,d^2\}\sm\{0\}))},\root{2h}\of{mN}\r]$$
 is at least $2^{q+h}\gs 2^{q'+h-1}$, and hence the interval $I_0=\{\al/2^{q'}:\al\in I\}$
 contains an integer $n_0$ such that $2^tm_0N_0-n_0^{2h}\not\in E$
 (by Lemma \ref{mn}). Thus $n=2^{q'}n_0\in I$ and
 $$mN-n^{2h}=4^{hq'}(2^tm_0N_0)-2^{2hq'}n_0^{2h}=4^{hq'}(2^tm_0N_0-n_0^{2h})\in\N\sm E.$$
 Applying Lemma \ref{Lem-SS} we see that there are $x,y,z,w\in\Z/N$ satisfying \eqref{mn^h}.
 As $n^{2h}\gs m(N-\min(\{a^2,b^2,c^2,d^2\}\sm\{0\}))$, we have $ax,by,cz,dw\gs0$ by Lemma \ref{Lem2.1}.
(If $a=0$ and $x<0$ we may use $-x>0$ to replace $x$
since $(-x)^2=x^2$. If $0\in\{b,c,d\}$ we also do such trivial replacements.)

In view of the above, we have completed the proof. \qed

\medskip
\noindent
{\it Proof of Theorem \ref{Th-abcd}}. If $r=0$, then it suffices to choose $x=y=z=w=0$.
If $r=m/n$ with $m,n\in\Z^+$, then by Lagrange's four-square theorem there are $x,y,z,w\in\N$
such that $mn=x^2+y^2+z^2+w^2$ and hence
$r=mn/n^2=(x/n)^2+(y/n)^2+(z/n)^2+(w/n)^2$. So Theorem \ref{Th-abcd} holds trivially if $r=0$
or $a=b=c=d=0$.

Below we assume $r>0$ and $a^2+b^2+c^2+d^2>0$.
Choose the least positive integer $q$
such that $a_0=aq,\, b_0=bq,\,c_0=cq,\,d_0=dq$ are all integers.
Write
$$\f r{q^2}=\f mn=\f{mn^{2h-1}}{n^{2h}}=\f{4^{hk}m_0}{n^{2h}},$$
where $m,m_0,n\in\Z^+$, $\gcd(m,n)=1$, $k\in\N$ and $4^h\nmid m_0$.
Choose any odd integer
$$M\gs\root{2h}\of{\f{M_h(a_0,b_0,c_0,d_0)}{m_0}}.$$
As $4^h\nmid m_0M^{2h}$ and $m_0M^{2h}\gs M_h(a_0,b_0,c_0,d_0)$, by Theorem \ref{Th-aux}
there are $x,y,z,w\in\N/(a_0^2+b_0^2+c_0^2+d_0^2)$ such that
$x^2+y^2+z^2+w^2=m_0M^{2h}$ and $a_0x+b_0y+c_0z+d_0w\in\Z^h$.
Let $\lambda =2^k/(Mn)$. Then
$$r=(\lambda^h q)^2(m_0M^{2h})=(\lambda^h qx)^2+(\lambda^h qy)^2
+(\lambda^h qz)^2+(\lambda^h qw)^2$$
and
$$a(\lambda^h qx)+b(\lambda^hqy)+c(\lambda^hqz)+d(\lambda^hqw)
=\lambda^h(a_0x+b_0y+c_0z+d_0w)\in\Q^h.$$
clearly
all the four rational numbers $\lambda^hqx,\lambda^hqy,\lambda^hqz,\lambda^hqw$
belong to the set
$$\f{q}{(a_0^2+b_0^2+c_0^2+d_0^2)(Mn)^h}\N
=\f{\N}{(a^2+b^2+c^2+d^2)q(Mn)^h}.$$
Note that $n\mid q^2$ if $r\in\Z^+$.

In view of the above, we have completed our proof of Theorem \ref{Th-abcd}. \qed

\section{Proofs of Theorems \ref{x+y-square}-\ref{four-square}}
 \setcounter{equation}{0}
 \setcounter{conjecture}{0}
 \setcounter{theorem}{0}

 \begin{lemma}\label{x+y=n^2} Let $h\in\{2,3\}$ and $m,n\in\N$ with $2m-n^{2h}\in \N\sm E$. Then
 there are $x,y,z,w\in\Z$ with $x^2+y^2+z^2+w^2=m$ and $x+y=n^h$.
 \end{lemma}
 \Proof. By the Gauss-Legendre theorem, we have $2m-n^{2h}=u^2+v^2+w^2$ for some $u,v,w\in\Z$.
 Clearly one of $u,v,w$ has the same parity with $n$. Without loss of generality, we assume that
 $w\eq n\pmod 2$ and hence $u\eq v\pmod 2$. Observe that
 $$m=\f{n^{2h}+w^2}2+\f{u^2+v^2}2=\l(\f{n^h+w}2\r)^2+\l(\f{n^h-w}2\r)^2+\l(\f{u+v}2\r)^2+\l(\f{u-v}2\r)^2$$
 with $(n^h+w)/2+(n^h-w)/2=n^h$. This proves the desired result. \qed

 \medskip
 \noindent{\it Proof of Theorem \ref{x+y-square}}. If $x,y,z\in\N/3$, $w\in\N$, and $x+y=t^2$ for some $t\in\Z$, then
 $16(x^2+y^2+z^2+w^2)=(4x)^2+(4y)^2+(4z)^2+(4w)^2$ with $4x,4y,4z\in\N/3$, $4w\in\N$ and $4x+4y=(2t)^2$.
 So it suffices to work with $m\in\N$ and $16\nmid m$.

 For $m\in\N$ with $16\nmid m$ and
 $$m\ls\l\lfloor\l(\f 4{\root4\of{2}-1}\r)^4\r\rfloor=199751,$$
 via a computer we find that $m$ can be written as $x^2+y^2+z^2+w^2$
 $(x,y,z,w\in\N)$ with $x+y$ a square except for $m=7,22,31,184,568,632$.
 Note that
 \begin{align*}7=&\l(\f13\r)^2+\l(\f 23\r)^2+\l(\f 73\r)^2+1^2\ \t{with}\ \f13+\f23=1^2;
 \\22=& \l(\f13\r)^2+\l(\f 23\r)^2+\l(\f 73\r)^2+4^2\ \t{with}\ \f13+\f23=1^2;
 \\31=& \l(\f13\r)^2+\l(\f 23\r)^2+\l(\f 73\r)^2+5^2\ \t{with}\ \f13+\f23=1^2;
 \\184=& \l(\f23\r)^2+\l(\f{10}3\r)^2+\l(\f {16}3\r)^2+12^2\ \t{with}\ \f23+\f{10}3=2^2;
 \\568=& \l(\f23\r)^2+\l(\f{10}3\r)^2+\l(\f{52}3\r)^2+16^2\ \t{with}\ \f23+\f{10}3=2^2;
 \\632=& \l(\f23\r)^2+\l(\f{10}3\r)^2+\l(\f{{20}}3\r)^2+24^2\ \t{with}\ \f23+\f{10}3=2^2.
 \end{align*}
 So any $m\in\N$ with $16\nmid m$ and $m\ls 199751$ has a desired representation.

 Below we suppose $m\in\N$, $16\nmid m$ and $m>199751$. Then
 $(\root 4\of 2-1)\root 4\of m\gs4$ and hence the interval $I=[\root 4\of{m},\root4\of{2m}]$
 contains four consecutive integers.

 As $32\nmid 2m$, by parts (i) and (ii) of Lemma \ref{mn}
 we can choose an integer $n\in I$ with $2m-n^4\not\in E$.

 Now that $2m-n^4\in\N\sm E$, by Lemma \ref{x+y=n^2} there are $x,y,z,w\in\Z$ such that
 $x^2+y^2+z^2+w^2=m$ and $x+y=n^2$. As $n^2\gs\sqrt{m(1^2+1^2-1^2)}$,
 applying Lemma \ref{Lem2.1} we get that $x\gs 0$ and $y\gs0$.

 In view of the above, we have completed the proof of Theorem \ref{x+y-square}. \qed

  \begin{lemma}\label{x+y+z=n^2} Let $m,n\in\N$ with $3m-n^4\in \N\sm E$. Then
 there are $x,y,z,w\in\Z$ with $x^2+y^2+z^2+w^2=m$ and $x+y+z=n^2$.
 \end{lemma}
 \Proof. Clearly $3m-n^4\not\eq1\pmod 3$. It is known that
 \begin{equation}\label{E236}\N\sm\{2x^2+3y^2+6z^2:\ x,y,z\in\Z\}=\{3q+1:\ q\in\N\}\cup E
 \end{equation}
 (cf. \cite[pp.\, 112-113]{D39}). Thus
 there are $w,u,v\in\Z$ such that
 $$3m-n^4=3w^2+6u^2+2(3v-n^2)^2=3(w^2+2u^2+2(3v^2-2n^2v))+2n^4.$$
 Hence
 $$m=w^2+2u^2+6v^2-4n^2v+n^4=w^2+x^2+y^2+z^2,$$
 where $x=u+v$, $y=v-u$ and $z=n^2-2z$. Note that $x+y+z=n^2$.
 This concludes the proof. \qed

 \medskip
 \noindent{\it Proof of Theorem \ref{x+y+z-square}}. It suffices to work with $m\in\N$ and $16\nmid m$. Via a computer we find that
 $m\in\N$ with $16\nmid m$ and
 $$m\ls\l\lfloor\l(\f2{\root4\of{3}-\root4\of{2}}\r)^4\r\rfloor=61762$$
 can be written as $x^2+y^2+z^2+w^2\ (x,y,z,w\in\N)$ with $x+y+z$ a square
except for $m$ in the set
\begin{align*}S=\{&3,\,13,\,18,\,21,\,23,\,40,\,47,\,56,\,68,\,73,\,75,\,79,\,83,\,125,\,149,
\\&157,\,180,\,193,\,205,\,312,\,328,\,411,\,424,\,431,\,501,\,708\}.
\end{align*}
For each $m\in S$, we find $x,y,z,w\in\N/27$ with $x^2+y^2+z^2+w^2=m$ such that $x+y+z$
is an integer square. For example,
\begin{align*}3=&\l(\f13\r)^2+\l(\f 7{27}\r)^2+\l(\f{11}{27}\r)^2+\l(\f{44}{27}\r)^2
\ \t{with}\ \f13+\f 7{27}+\f{11}{27}=1^2;
\\13=&\l(\f13\r)^2+\l(\f 2{9}\r)^2+\l(\f{4}{9}\r)^2+\l(\f{32}{9}\r)^2
\ \t{with}\ \f13+\f 2{9}+\f{4}{9}=1^2;
\\18=&2^2+\l(\f13\r)^2+\l(\f 53\r)^2+\l(\f{10}3\r)^2\ \t{with}\ 2+\f13+\f53=2^2;
\\21=&0^2+\l(\f 8{27}\r)^2+\l(\f{19}{27}\r)^2+\l(\f{122}{27}\r)^2
\ \t{with}\ 0+\f 8{27}+\f{19}{27}=1^2;
\\708=&1^2+\l(\f 53\r)^2+\l(\f{67}3\r)^2+\l(\f{43}3\r)^2
\ \t{with}\ 1+\f 53+\f{67}3=5^2.
\end{align*}

 Now we handle the case $m>61762$. As $(\root 4\of3-\root4\of2)\root 4\of{m}\gs2$,
 the interval $I=[\root4\of{2m},\root4\of{3m}]$ contains two consecutive integers.
 Choose an integer $n\in I$ such that
 $$n\ \t{is odd}\iff m\eq2\pmod4\ \t{or}\ m\eq 4,5\pmod 8.$$
 Then $3m-n^4\in\N\sm E$ by Lemma \ref{mn}(i). Thus, by Lemma \ref{x+y+z=n^2}
  there are $x,y,z,w\in\Z$
 such that $3m-n^4=x^2+y^2+z^2+w^2$ and $x+y+z=n^2$. As $n^2\ge\sqrt{(1^2+1^2+1^2-1^2)m}$,
 we have $x,y,z\in\N$ by applying Lemma \ref{Lem2.1}.
 This concludes our proof. \qed

\medskip
\noindent
{\it Proof of Theorem \ref{four-square}}. It suffices to work with $m\in\N$ and $16\nmid m$.
For $m\in\N$ with $16\nmid m$ and
$$m\ls\l\lfloor\l(\f4{\root4\of{4}-\root4\of{3}}\r)^4\r\rfloor=2759711,$$
via a computer we find that $m$ can be written as $x^2+y^2+z^2+w^2$
with $x,y,z,w\in\N/50$ such that $x+y+z+w$ is a rational square.
(Actually we may even require $x,y,z,w\in\N$ if $57436<m\ls 2759711$.)
For example,
$$2=\l(\f 45\r)^2+\l(\f 45\r)^2+\l(\f3{25}\r)^2+\l(\f{21}{25}\r)^2
\ \t{with}\ \f45+\f45+\f3{25}+\f{21}{25}=\l(\f 85\r)^2,$$
\begin{align*}11=&\l(\f1{50}\r)^2+\l(\f7{50}\r)^2+\l(\f{111}{50}\r)^2+\l(\f{123}{50}\r)^2,
\\& \f1{50}+\f7{50}+\f{111}{50}+\f{123}{50}=\l(\f{11}5\r)^2;
\\24992=&0^2+\l(\f{752}5\r)^2+\l(\f{228}{25}\r)^2+\l(\f{1196}{25}\r)^2,
\\&0+\f{752}{5}+\f{228}{25}+\f{1196}{25}=\l(\f{72}5\r)^2;
\\56956=&\l(\f 3{25}\r)^2+\l(\f{1561}{25}\r)^2+\l(\f{1803}{25}\r)^2+\l(\f{5469}{25}\r)^2,
\\&\f3{25}+\f{1561}{25}+\f{1803}{25}+\f{5469}{25};
\\57436=&\l(\f1{25}\r)^2+\l(\f{19}{25}\r)^2+\l(\f{1603}{25}\r)^2+\l(\f{5773}{25}\r)^2,
\\&\f1{25}+\f{19}{25}+\f{1603}{25}+\f{5773}{25}=\l(\f{86}5\r)^2.
\end{align*}

Now we handle the case $m>2759711$. As $(\root4\of4-\root4\of3)\root4\of m\gs4$,
the interval $I=[\root4\of{3m},\root4\of{4m}]$ contains four consecutive integers.
As $64\nmid 4m$, by Lemma \ref{mn} we may choose $n\in I$ such that $4m-n^4\in\N\sm E$.
Write $4m-n^4=x_1^2+x_2^2+x_3^2$ with $x_1,x_2,x_3\in\Z$.
If $m$ is odd, then $4m-n^4\eq 3\pmod 8$ and hence $x_1,x_2,x_3\eq1\eq n^2\pmod 2$.
If $m$ is even, then $4m-n^4\eq0\pmod 8$ and hence $x_1,x_2,x_3\eq 0\eq n^2\pmod 2$.
So there are $u,v,w\in\Z$ such that
\begin{equation}\label{uvw}
4m-n^4=(2u-n^2)^2+(2v-n^2)^2+(2w-n^2)^2.\end{equation}
 If $m$ is even, then \eqref{uvw} modulo $8$ yields
$0\eq 4(u^2+v^2+w^2)\pmod 8$ and hence $u+v+w\eq0\pmod2$.
When $m$ and $u+v+w$ are both odd, for $u'=n^2-u$ we have $(2u'-n^2)^2=(2u-n^2)^2$
and $u'+v+w\eq0\pmod2$. So, without loss of generality, we may assume that $u+v+w\eq0\pmod2$
no matter $m$ is even or not. Set
$$x=\f{-u+v+w}2,\ y=\f{u-v+w}2,\ z=\f{u+v-w}2.$$
Then $x+y=w$, $y+z=u$ and $x+z=v$. Hence
\begin{align*}4m=&(2(y+z)-n^2)^2+(2(x+z)-n^2)^2+(2(x+y)-n^2)^2+n^4
\\=&4(x^2+y^2+z^2+(x+y+z-n^2)^2)
\end{align*}
and thus $m=x^2+y^2+z^2+t^2$ with $t=n^2-(x+y+z)$.
As $n^2\gs\sqrt{(1^2+1^2+1^2+1^2-1^2)m}$, by Lemma \ref{Lem2.1} we have $x,y,z,t\in\N$.

In view of the above, we have completed the proof of Theorem \ref{four-square}. \qed

\section{Proofs of Theorems \ref{Th-112}-\ref{Th122}}
\setcounter{equation}{0}
\setcounter{conjecture}{0}

 \begin{lemma}\label{Lem-112} Let $m,n\in\Z^+$ with $6m-n^4\in\N\sm E$. Then
 $m$ can be written as $x^2+y^2+z^2+w^2$ with $x,y,z,w\in\Z$ and $x+y+2z=n^2$.
 \end{lemma}
 \Proof. Clearly $6m-n^4\not\eq 1\pmod3$. In view of \eqref{E236}, there are $u,v,w\in\Z$
 such that $6m-n^4=6w^2+3(2u+n^2)^2+2(3v+n^2)^2$. (Note that $(3v-n^2)^2=(3(-v)+n^2)^2$.)
 It follows that
 $$m=w^2+2u^2+2n^2(u+v)+n^4=w^2+(y+z+n^2)^2+(z-y)^2+(-z)^2$$
 with $(y+z+n^2)+(z-y)+2(-z)=n^2$. This concludes the proof. \qed

 \medskip
\noindent{\it Proof of Theorem \ref{Th-112}}. It suffices to
work with $m\in\N$ with $16\nmid m$.

Via a computer we find that $m\in\N$ with $16\nmid m$ and
$$m\ls \l\lfloor\l(\f4{\root4\of{6}-\root4\of{5}}\r)^4\r\rfloor=10824724$$
can be written as $x^2+y^2+z^2+w^2\ (x,y,z,w\in\N)$ with $x+y+2z\in\Z^2$
except for $m=58,\,168,\,1576,\,1640$. Note that
\begin{align*}58=&1^2+\l(\f{5}3\r)^2+\l(\f23\r)^2+\l(\f{22}3\r)^2
\ \t{with}\ 1+\f{5}3+2\times\f23=2^2;
\\168=&0^2+\l(\f 83\r)^2+\l(\f 23\r)^2+\l(\f{38}3\r)^2
\ \t{with}\ 0+\f 83+2\times\f23=2^2;
\\1576=&0^2+\l(\f{64}3\r)^2+\l(\f{22}3\r)^2+\l(\f{98}3\r)^2
\ \t{with}\ 0+\f{64}3+2\times\f{22}3=6^2;
\\1640=&6^2+\l(\f{58}3\r)^2+\l(\f{16}3\r)^2+\l(\f{104}3\r)^2
\ \t{with}\ 6+\f{58}3+2\times\f{16}3=6^2.
\end{align*}

Now we handle the case $m>10824724$. As $(\root4\of{4}-\root4\of{3})\root4\of m\gs4$,
the length of interval $I=\root4\of{3m},\root4\of{4m}$ is at least four.
As $32\nmid 6m$, by Lemma \ref{mn} there is an integer $n\in I$ such that $6m-n^4\not\in E$.
In view of Lemma \ref{Lem-112}, there are $x,y,z,w\in\Z$ such that
$x^2+y^2+z^2+w^2=m$ and $x+y+2z=n^2$.  As $n^2\gs\sqrt{(1^2+1^2+2^2-1^2)m}$, by Lemma \ref{Lem2.1} we have $x,y,z\in\N$.

In view of the above, we have completed the proof of Theorem \ref{Th-112}. \qed

 \begin{lemma}\label{Lem2.2} Let $m,n\in\Z^+$ with $3\nmid n$ and $9m-n^4\in\N\sm E$. Then
 $m$ can be written as $x^2+y^2+z^2+w^2$ with $x,y,z,w\in\Z$ and $x+2y+2z=n^2$.
 \end{lemma}
 \Proof. By the Gauss-Legendre theorem, $9m-n^4=a^2+b^2+c^2$ for some $a,b,c\in\Z$.
 As $3\nmid n$, we have $3\mid abc$. Without loss of generality, we suppose $c=3w$ with $w\in\Z$.
 Since $a^2+b^2\eq-n^4\eq2\pmod 3$, we have $3\nmid ab$. Without loss of generality, we may assume that
 $a\eq2\pmod 3$ and $b\eq-2\pmod 3$ (otherwise we may change the signs of $a$ and $b$ suitably).
 Clearly, $a=3u+2n^2$ and $b=3v-2n^2$ for some $u,v\in\Z$. Observe that
 $$12n^2(u-v)+8n^4\eq (3u+2n^2)^2+(3v-2n^2)^2=a^2+b^2\eq -n^4\pmod 9$$
 and hence $u\eq v\pmod 3$ since $3\nmid n$. Set
 $$y=-\f{2u+v}3\ \ \t{and}\ \ z=\f{u+2v}3.$$
 Then
 \begin{align*}9m-n^4=&a^2+b^2+c^2=(3u+2n^2)^2+(3v-2n^2)^2+9w^2
 \\=&(3(-2y-z)+2n^2)^2+(3(2z+y)-2n^2)^2+9w^2
 \\=&9(2y+z)^2+9(2z+y)^2-36n^2(y+z)+8n^4+9w^2
 \end{align*}
 and hence
 $$m=n^4+w^2+(2y+z)^2+(2z+y)^2-4n^2(y+z)=x^2+y^2+z^2+w^2,$$
 where $x=n^2-2(y+z)$. Note that $x+2y+2z=n^2$ as desired. \qed

\medskip
\noindent{\it Proof of Theorem \ref{Th122}}. It suffices to
work with $m\in\N$ with $16\nmid m$.

Via a computer we find that $m\in\N$ with $16\nmid m$ and
$$m\ls \l\lfloor\l(\f2{\root4\of{9}-\root4\of{8}}\r)^4\r\rfloor=40125453$$
can be written as $x^2+y^2+z^2+w^2$ $(x,y,z,w\in\N)$ with $x+2y+2z$ a square
except for $m=7,24,55,112,120,255$. Note that
\begin{align*}7=&\l(\f1{13}\r)^2+\l(\f1{13}\r)^2+\l(\f5{13}\r)^2+\l(\f{34}{13}\r)^2
,\  \f1{13}+2\times\f1{13}+2\times\f 5{13}=1^2;
\\24=&\l(\f{24}{13}\r)^2+\l(\f4{13}\r)^2+\l(\f{10}{13}\r)^2+\l(\f{58}{13}\r)^2
,\  \f{24}{13}+2\times\f4{13}+2\times\f {10}{13}=2^2;
\\55=&\l(\f{11}{13}\r)^2+\l(\f{10}{13}\r)^2+\l(\f{43}{13}\r)^2+\l(\f{85}{13}\r)^2
,\  \f{11}{13}+2\times\f{10}{13}+2\times\f {43}{13}=3^2;
\\120=&\l(\f{12}{13}\r)^2+\l(\f{34}{13}\r)^2+\l(\f{64}{13}\r)^2+\l(\f{122}{13}\r)^2
,\  \f{12}{13}+2\times\f{34}{13}+2\times\f {64}{13}=4^2;
\\255=&\l(\f{3}{13}\r)^2+\l(\f{7}{13}\r)^2+\l(\f{154}{13}\r)^2+\l(\f{139}{13}\r)^2
,\  \f{3}{13}+2\times\f{7}{13}+2\times\f {154}{13}=5^2.
\end{align*}

Now we handle the case $m>40125453$. As $(\root4\of{9}-\root4\of{8})\root4\of{m}\gs4$,
the interval $I=[\root 4\of{8m},\root4\of{9m}]$ contains four consecutive integers.
Note that for any $s\in\Z$ either $s$ or $s+2$ is not divisible by $3$.
We may choose an integer $n\in I$ with $3\nmid n$ such that $n$ is odd if and only if $m\eq2\pmod4$
or $m\eq4,7\pmod 8$. Then $9m-n^4\in\N\sm E$ by Lemma \ref{mn}.
In light of Lemma \ref{Lem2.2}, there are $x,y,z,w\in\Z$ such that $x^2+y^2+z^2+w^2=m$ and
$x+2y+2z=n^2$. As $n^2\gs\sqrt{(1^2+2^2+2^2-1^2)m}$, by Lemma \ref{Lem2.1} we have $x,y,z\in\N$.

In view of the above, we have completed the proof of Theorem \ref{Th122}. \qed

\section{Proofs of Theorems \ref{Th123}-\ref{Th-1123}}
\setcounter{equation}{0}
\setcounter{conjecture}{0}

 \begin{lemma}\label{Lem2.3} Let $m,n\in\Z^+$.

 {\rm (i)} If $11m-n^4\in\N\sm E$, then there are $x,y,z,w\in\Z$
 such that $x^2+y^2+z^2+w^2=m$ and $x+y+3z=n^2$.
 If $14m-n^4\in\N\sm E$, then there are $x,y,z,w\in\Z$
 such that $x^2+y^2+z^2+w^2=m$ and $x+2y+3z=n^2$.

 {\rm (ii)} If $7m-n^4\in\N\sm E$, then there are $x,y,z,w\in\Z$
 such that $x^2+y^2+z^2+w^2=m$ and $x+y+z+2w=n^2$.
 If $15m-n^4\in\N\sm E$, then there are $x,y,z,w\in\Z$
 such that $x^2+y^2+z^2+w^2=m$ and $x+y+2z+3w=n^2$.
 If $23m-n^4\in\N\sm E$, then there are $x,y,z,w\in\Z$
 such that $x^2+y^2+z^2+w^2=m$ and $x+y+2z+3w=n^2$.
 \end{lemma}
 \Proof. For the five Hamilton quaternions
 \begin{gather*}\zeta_1=1+i+3j+0k,\ \zeta_2=1+2i+3j+0k,\ \zeta_3=1+i+j+2k,
 \\\zeta_4=1+i+2j+3k,\ \zeta_5=1+2i+3j+3k
 \end{gather*}
 their norms are $11,\,14,\,7,\,15,\,23$ respectively. Note that each of $11,\,14,\,7,\,15,\,23$
 can be written uniquely in the form $a^2+b^2+c^2+d^2$ with $a,b,c,d\in\N$ and $a\ls b\ls c\ls d$.
 Applying \cite[Theorem 2]{MT}, we immediately get the desired result. \qed

\medskip
\noindent{\it Proof of Theorem \ref{Th123}}. It suffices to work with $m\in\N$ and $16\nmid m$.

(i) Via a computer, $m\in\N$ with $16\nmid m$ and
$$m\ls\l\lfloor\l(\f 2{\root 4\of{11}-\root4\of{10}}\r)^4\r\rfloor =4732224$$
can be written as $x^2+y^2+z^2+w^2\ (x,y,z,w\in\N)$ with $x+y+3z$ a square, except for
$m=7,\,15,\,133$. Note that
\begin{align*}7=&\l(\f23\r)^2+\l(\f 73\r)^2+\l(\f 13\r)^2+1^2
\ \t{with}\ \f23+\f73+3\times\f13=2^2,
\\15=&\l(\f 23\r)^2+\l(\f 73\r)^2+\l(\f 13\r)^2+3^2
\ \t{with}\ \f 23+\f73+3\times\f13=2^2,
\\133=&\l(\f 23\r)^2+\l(\f{16}3\r)^2+\l(\f{19}3\r)^2+8^2\ \t{with}\ \f 23+\f{16}3+3\times\f{19}3=5^2.
\end{align*}

Now we handle the case $m>4732224$. As $(\root4\of{11}-\root4\of{10})\root4\of{m}\gs2$,
the interval $I=[\root4\of{10m},\root4\of{11m}]$ contains two consecutive integers.
As $16\nmid 11m$, by Lemma \ref{mn} we can choose an integer $n\in I$ such that $11m-n^4\in \N\sm E$.
In view of Lemma \ref{Lem2.3}, there are $x,y,z,w\in\N$ such that $x^2+y^2+z^2+w^2=m$
and $x+y+3z=n^2$. As $n^2\gs\sqrt{(1^2+1^2+3^2-1^2)m}$, by Lemma \ref{Lem2.1} we have
$x,y,z\in\N$.

(ii) Via a computer, $m\in\N$ with $16\nmid m$ and
$$m\ls\l\lfloor\l(\f 4{\root 4\of{14}-\root4\of{13}}\r)^4\r\rfloor = 161049608$$
can be written as $x^2+y^2+z^2+w^2\ (x,y,z,w\in\N)$ with $x+2y+3z$ a square, except for
$m=12,\,72,\,76,\,92$. Note that
\begin{align*}12=&\l(\f25\r)^2+\l(\f 65\r)^2+\l(\f 25\r)^2+\l(\f{16}5\r)^2
\ \t{with}\ \f25+2\times\f 65+3\times25=2^2,
\\72=&\l(\f 45\r)^2+\l(\f 25\r)^2+\l(\f 45\r)^2+\l(\f{42}5\r)^2
\ \t{with}\ \f 45+2\times25+3\times45=2^2,
\\76=&\l(\f 75\r)^2+7^2+\l(\f15\r)^2+5^2\ \t{with}\ \f 75+2\times7+3\times\f15=4^2.
\end{align*}

Now we handle the case $m>161049608$. As $(\root4\of{14}-\root4\of{13})\root4\of{m}\gs4$,
the interval $I=[\root4\of{13m},\root4\of{14m}]$ contains four consecutive integers.
As $32\nmid 14m$, by Lemma \ref{mn} we can choose an integer $n\in I$ such that $14m-n^4\in \N\sm E$.
In view of Lemma \ref{Lem2.3}, there are $x,y,z,w\in\N$ such that $x^2+y^2+z^2+w^2=m$
and $x+2y+3z=n^2$. As $n^2\gs\sqrt{(1^2+2^2+3^2-1^2)m}$, by Lemma \ref{Lem2.1} we have
$x,y,z\in\N$.

By the above, we have finished the proof of Theorem \ref{Th123}. \qed

\medskip
\noindent{\it Proof of Theorem \ref{Th1.2}}. In view of Lemma \ref{Lem2.3}, it suffices to find $a\in\Z^+$
with $14m-2^{4a}\in\N\sm E$. If $14m-2^{4a}\in\N\sm E$ then $14(16m)-2^{4(a+1)}=4^2(14m-2^{4a})\in\N\sm E$. Note also that
$$16=4^2+0^2+0^2+0^2\ \t{and}\ 160=4^2+0^2+0^2+12^2\ \t{with}\ 4+2\times0+3\times0=4^1.$$
So we only need to handle the case $m\not\eq0\pmod {16}$.

If $m\in\{2,\ldots,18\}\sm\{10,16\}$, then we can verify the desired result directly.

Below we assume that $m>18$ and $16\nmid m$. Note that $14m\gs14\times 19=266>2^8$.

If $2\nmid m$, then $14m-2^4\in\N\sm E$ since $14m\eq2\pmod 4$.

In the case $m\eq2\pmod4$, we have $14m=4q$ for some odd integer $q>64$.
If $q\not\eq 7\pmod 8$, then $14m-2^8=4(q-64)\in\N\sm E$.
If $q\eq7\pmod 8$, then $14m-2^4=4(q-4)\in\N\sm E$.

If $m\eq4\pmod 8$, then $14m\eq 8\pmod {16}$ and hence $14m-2^4\in\N\sm E$.

In the case $m\eq 8\pmod{16}$, we have $14m=16q$ for some odd integer $q>16$.
If $q\not\eq 7\pmod 8$, then $14m-2^8=16(q-16)\in \N\sm E$.
If $q\eq7\pmod 8$, then $14m-2^4=16(q-1)\in \N\sm E$.

In view of the above, we have completed the proof of Theorem \ref{Th1.2}. \qed

\medskip
\noindent{\it Proof of Theorem \ref{Th-1123}}. It suffices to work with $m\in\N$ and $16\nmid m$.

(i) Via a computer we find that $m\in\N$ with $16\nmid m$ and
$$m\ls\l\lfloor\l(\f2{\root4\of{7}-\root4\of{6}}\r)^4\r\rfloor=1119041$$
can be written as $x^2+y^2+z^2+w^2\ (x,y,z,w\in\N)$ with $x+y+z+2w\in\Z^2$
except for values of $m$ among
\begin{gather*}2,\,7,\,9,\,11,\,19,\,23,\,24,\,25,\,31,\,34,\,36,\,49,
\,57,\,62,\,68,\,72,\,73,\,74,\,82,\,89,
\\116,\,119,\,135,\,139,\,143,\,166,\,167,
\,168,\,179,\,184,\,268,\,292,\,
\\340,\,424,\,472,\,552,\,568,\,583,\,863,\,2744.
\end{gather*}
All these numbers can be written as $x^2+y^2+z^2+w^2$
with $x,y,z\in\N/7\cup \N/9\cup \N/25$ and $w\in\N$ such that $x+y+z+2w\in\Q^2$.
For example,
\begin{align*}7=&\l(\f 75\r)^2+\l(\f{11}{25}\r)^2+\l(\f{23}{25}\r)^2=2^2,
\ \f 75+\f{11}{25}+\f{23}{25}+2\times2=\l(\f{13}5\r)^2;
\\340=&\l(\f 57\r)^2+\l(\f{31}7\r)^2+\l(\f{125}7\r)^2+1^2,
\ \f57+\f{31}7+\f{125}7+2\times1=5^2;
\\863=&\l(\f 53\r)^2+\l(\f{23}3\r)^2+\l(\f{83}3\r)^2+6^2,
\ \f53+\f{23}3+\f{83}3+2\times6=7^2;
\\2744=&\l(\f29\r)^2+\l(\f{160}9\r)^2+\l(\f{442}9\r)^2+4^2,
\ \f29+\f{160}9+\f{442}9+2\times4=\l(\f{26}3\r)^2.
\end{align*}

Now we handle the case $m>1119041$. As $(\root4\of{7}-\root4\of{6})\root4\of m\gs2$,
the interval $I=[\root4\of{6m},\root4\of{7m}]$ contains two consecutive integers.
As $16\nmid 7m$, by Lemma \ref{mn} there is an integer $n\in I$ such that $7m-n^4\in\N\sm E$.
In light of Lemma \ref{Lem2.3}(ii), there are $x,y,z,w\in\Z$ such that
$x^2+y^2+z^2+w^2=m$ and $x+y+z+2w=n^2$. Since $n^2\gs\sqrt{(1^2+1^2+1^2+2^2-1^2)m}$,
applying Lemma \ref{Lem2.1} we find that $x,y,z,w\in\N$.

(ii) Via a computer we find that $m\in\N$ with $16\nmid m$ and
$$m\ls\l\lfloor\l(\f2{\root4\of{15}-\root4\of{14}}\r)^4\r\rfloor=12474176$$
can be written as $x^2+y^2+z^2+w^2\ (x,y,z,w\in\N)$ with $x+y+2z+3w\in\Z^2$
except for $m=11,\,15,\,25,\,33,\,36,\,40,\,71,\,79,\,97,\,127,\,153.$
Note that
\begin{align*}11=&\l(\f15\r)^2+\l(\f35\r)^2+\l(\f{16}5\r)^2+\l(\f 35\r)^2
,\  \f15+\f35+2\times\f{16}5+3\times\f 35=3^2;
\\15=&3^2+\l(\f 75\r)^2+2^2+\l(\f15\r)^2
,\  3+\f 75+2\times2+3\times\f15=3^2;
\\25=&0^2+\l(\f{24}5\r)^2+0^2+\l(\f 75\r)^2,\ 0+\f{24}5+2\times0+3\times\f 75=3^2;
\\
33=&\l(\f 25\r)^2+\l(\f{28}5\r)^2+\l(\f 65\r)^2+\l(\f15\r)^2
,\ \f25+\f{28}5+2\times\f 65+3\times\f15=3^2;
\\
36=&\l(\f 75\r)^2+\l(\f{29}5\r)^2+\l(\f35\r)^2+\l(\f15\r)^2,
\ \f 75+\f{29}5+2\times\f35+3\times\f15=3^2;
\end{align*}
\begin{align*}
40=&0^2+\l(\f{26}5\r)^2+0^2+\l(\f{18}5\r)^2,
\ 0+\f{26}5+2\times0+3\times\f{18}5=4^2;
\\71=&1^2+\l(\f{27}5\r)^2+6^2+\l(\f{11}5\r)^2,
\ 1+\f{27}5+2\times6+3\times\f{11}5=5^2;
\\79=&\l(\f25\r)^2+\l(\f{41}5\r)^2+\l(\f{17}5\r)^2+\l(\f15\r)^2,
\ \f25+\f{41}5+2\times\f{17}5+3\times\f15=4^2;
\\
97=&0^2+\l(\f{38}5\r)^2+6^2+\l(\f 95\r)^2,\ 0+\f{38}5+2\times6+3\times\f 95=5^2;
\\127=&\l(\f35\r)^2+\l(\f 65\r)^2+\l(\f{51}5\r)^2+\l(\f{23}5\r)^2,
\ \f35+\f65+2\times\f{51}5+3\times\f{23}5=6^2;
\\153=&\l(\f 25\r)^2+\l(\f{59}5\r)^2+\l(\f{14}5\r)^2+\l(\f{12}5\r)^2,
\ \f25+\f{59}5+2\times\f{14}5+3\times\f{12}5=5^2.
\end{align*}

Now we handle the case $m>12474176$. As $(\root4\of{15}-\root4\of{14})\root4\of m\gs2$,
the interval $I=[\root4\of{14m},\root4\of{15m}]$ contains two consecutive integers.
As $16\nmid 15m$, by Lemma \ref{mn} there is an integer $n\in I$ such that $15m-n^4\in\N\sm E$.
In light of Lemma \ref{Lem2.3}(ii), there are $x,y,z,w\in\Z$ such that
$x^2+y^2+z^2+w^2=m$ and $x+y+2z+3w=n^2$. Since $n^2\gs\sqrt{(1^2+1^2+2^2+3^2-1^2)m}$,
applying Lemma \ref{Lem2.1} we find that $x,y,z,w\in\N$.

(iii) Via a computer we find that $m\in\N$ with $16\nmid m$ and
$$m\ls\l\lfloor\l(\f2{\root4\of{23}-\root4\of{22}}\r)^4\r\rfloor=46645840$$
can be written as $x^2+y^2+z^2+w^2\ (x,y,z,w\in\N)$ with $x+2y+3z+3w\in\Z^2$
except for values of $m$ among
\begin{gather*}3,\,7,\,8,\,23,\,25,\,40,\,43,\,55,\,56,\,97,\,120,\,168,\,172,\,376,
\end{gather*}
All these numbers can be written as $x^2+y^2+z^2+w^2$ with $y\in\N$
and $x,z,w\in\N/5\cup\N/7\cup\N/17$ such that $x+2y+3z+3w\in\Z^2$. For example,
\begin{align*}
3=&\l(\f75\r)^2+1^2+0^2+\l(\f15\r)^2,\ \f 75+2\times1+3\times0+3\times\f15=2^2;
\\55=&\l(\f{50}7\r)^2+1^2+\l(\f 57\r)^2+\l(\f{11}7\r)^2
,\ \f{50}7+2\times1+3\times\f 57+3\times\f{11}7=4^2;
\\56=&\l(\f{122}{17}\r)^2+0^2+\l(\f{20}{17}\r)^2+\l(\f{30}{17}\r)^2
,\  \f{122}{17}+2\times0+3\times\f{20}{17}+3\times\f{30}{17}=4^2;
\\376=&\l(\f{64}5\r)^2+4^2+\l(\f25\r)^2+14^2,\ \f{64}5+2\times4+3\times\f25+3\times14=8^2.
\end{align*}

Now we handle the case $m>46635840$. As $(\root4\of{23}-\root4\of{22})\root4\of m\gs2$,
the interval $I=[\root4\of{22m},\root4\of{23m}]$ contains two consecutive integers.
As $16\nmid 23m$, by Lemma \ref{mn} there is an integer $n\in I$ such that $23m-n^4\in\N\sm E$.
In light of Lemma \ref{Lem2.3}(ii), there are $x,y,z,w\in\Z$ such that
$x^2+y^2+z^2+w^2=m$ and $x+2y+3z+3w=n^2$. Since $n^2\gs\sqrt{(1^2+2^2+3^2+3^2-1^2)m}$,
applying Lemma \ref{Lem2.1} we find that $x,y,z,w\in\N$.

In view of the above, we have completed the proof of Theorem \ref{Th-1123}. \qed

\section{Some conjectures}
\setcounter{equation}{0}
\setcounter{conjecture}{0}

\begin{conjecture} [4-4-2-2 Conjecture] Let $a,b\in\Z^+$ with $\gcd(a,b)=1$. Then each $r\in\Q_{\gs0}$ can be written as
$aw^4+bx^4+y^2+z^2$ with $w,x,y,z\in\Q$. Moreover, for any $n\in\Z^+$
there is a positive integer $m<2(a+b)$ such that $m^4n=2^{4w}+x^4+y^2+z^2$ for some $w,x,y,z\in\N$.
\end{conjecture}
\begin{remark} For example, we conjecture that for any $n\in\Z^+$ we can write $3^4n$
as $2^{4w}+x^4+y^2+z^2$ with $w,x,y,z\in\N$; this has been verified for $n\ls 10^8$.
Moreover, for $a,b\in\{1,\ldots,10\}$ with $a+b>2$ and $\gcd(a,b)=1$
we conjecture that for any $n\in\Z^+$ we can write $m^4n=a2^{4w}+bx^4+y^2+z^2$ with $w,x,y,z\in\N$,
provided that $(a,b,m)$ is among the following ordered triples
\begin{align*}&(1,2,4),(1,3,4),(1,4,7),(1,5,8),(1,6,3),(1,7,4),(1,8,6),(1,9,9),(1,10,6),
\\&(2,1,4),(2,3,3),(2,5,3),(2,7,6),(2,9,6),(3,1,4),(3,2,4),(3,4,4),(3,5,4),
\\&(3,7,4),(3,8,2),(3,10,6),(4,1,7),(4,3,4),(4,5,17),(4,7,4),(4,9,9),(5,1,8),
\\&(5,2,3),(5,3,6),(5,4,11),(5,6,6),(5,7,4),(5,8,6),(5,9,6),(6,1,4),(6,5,2),
\\&(6,7,4),(7,1,4),(7,2,4),(7,3,4),(7,4,4),(7,5,4),(7,6,6),(7,8,4),(7,9,6),
\\&(7,10,6),(8,1,8),(8,3,4),(8,5,12),(8,7,2),(8,9,12),(9,1,3),(9,2,3),(9,4,9),
\\&(9,5,6),(9,7,6),(9,8,6),(9,10,3),(10,1,5),(10,3,6),(10,7,6),(10,9,6);
\end{align*}
we have verified this for $n\ls2\times10^5$. It seems that the least positive integer $m$
such that for any $n\in\Z^+$ we can write $m^4n=16\times2^{4w}+33x^4+y^2+z^2$ with $w,x,y,z\in\N$,
is $2\times33=66$.
\end{remark}

\begin{conjecture} Let $a,b\in\Z^+$ with $\gcd(a,b)=1$ and $\{a\}_3+\{b\}_3\not=1$,
 where $\{c\}_3$ denotes the least nonnegative residue of an integer $c$ modulo $3$.
 Then each $r\in\Q_{\gs0}$ can be written as
$aw^4+bx^4+y^2+3z^2$ with $w,x,y,z\in\Q$. Moreover, for any $n\in\Z^+$
there is a positive integer $m\ls2(a+b+1)$ such that $m^4n=2^{4w}+x^4+y^2+3z^2$ for some $w,x,y,z\in\N$.
\end{conjecture}
\begin{remark} For example, we conjecture that for any $n\in\Z^+$ we can write $5^4n$
as $2^{4w}+x^4+y^2+3z^2$ with $w,x,y,z\in\N$, and also we can write $8^4n$ as
$2^{4w}+2x^4+y^2+3z^2$ with $w,x,y,z\in\N$.
\end{remark}

\begin{conjecture}  Let $K$ be a totally real field (i.e,, an algebraic subfield of the field $\R$ of real numbers), and set $K_{\gs0}=\{t\in K:\ t\gs0\}$.

{\rm (i)} Let
$a,b,c,d\in K$ with $a>0$, and let $h\in\{2,3\}$. Then each $r\in K_{\gs0}$ can be written
as $x^2+y^2+z^2+w^2$ with $x,y,z,w\in K_{\gs0}$ such that $ax+by+cz+dw\in\{t^h:\ t\in K\}$.

{\rm (ii)} Each $r\in K_{\gs0}$ can be written as $x^4+y^4+z^2+w^2$ with $x,y,z,w\in K$.
\end{conjecture}
\begin{remark} In 1921 C.L. Siegel \cite{Si}
confirmed a conjecture of D. Hilbert by proving that each nonnegative element of a totally real field $K$ can be written as $x^2+y^2+z^2+w^2$ with $x,y,z,w\in K$. Parts (i) and (ii) are motivated by
Theorem 1.1 and the 4-4-2-2 Conjecture respectively.
\end{remark}

\begin{conjecture} [1-2-3-4 Conjecture] Any $n\in\N$ with $n\not=158$ can be written as
$w^2+2x^2+3y^4+4z^4$ with $w,x,y,z\in\N$.
\end{conjecture}
\begin{remark} This has been verified for $n\ls 10^8$. See \cite[A346643]{OEIS} for related data. For example, $4254$ has a unique required
representation: $4254=45^2+2\times31^2+3\times3^4+4\times2^4$.
We also conjecture that $744$ is the only nonnegative integer which cannot be written as 
$w^2+2x^2+y^4+4z^4$ with $w,x,y,z\in\N$ (cf. \cite[A347]{OEIS}).
\end{remark}

\begin{conjecture} Each $n\in\N$ with $n\not=95,255$ can be written as $w^2+x^2+y^4+2z^4$
$(w,x,y,z\in\N)$ with $y$ or $z$ an integer square.
\end{conjecture}
\begin{remark} See \cite[A347857]{OEIS} for related data.
We also conjecture that any integer $n>8640$ can be written as $w^2+4x^2+y^4+2z^6$
with $w,x,y,z\in\N$, and that any integer $n>20319$ can be written as $w^2+4x^2+y^4+2z^8$
with $w,x,y,z\in\N$.
\end{remark}

\begin{conjecture} Let $a$ be a positive integer.

{\rm (i)} If $a\eq2\pmod4$, then any sufficiently large integer can be written as $aw^4+x^4+(2y)^2+z^2$
with $w,x,y,z\in\N$.

{\rm (ii)} If $a\eq1\pmod2$, then any sufficiently large integer can be written as $aw^4+2x^4+(2y)^2+z^2$
with $w,x,y,z\in\N$.
\end{conjecture} 
\begin{remark} For $a\eq2\pmod4$ let $M(a)$ be the largest integer not of the form
$aw^4+x^4+(2y)^2+z^2\ (w,x,y,z\in\N)$. Based on our computation, we guess that 
\begin{gather*}M(2)=255,\ M(6)=2716,\ M(10)=598,\ M(14)=8427,
\\ M(18)= 2463,\ M(22)=3884, M(26)=14988,\ M(30)=10843.
\end{gather*}
For $a\eq1\pmod2$ let $N(a)$ be the largest integer not of the form
$aw^4+2x^4+(2y)^2+z^2\ (w,x,y,z\in\N)$. We conjecture that
\begin{gather*}N(1)=255,\ N(3)=303,\ N(5)=497,\ N(7)=3182,\ N(9)= 4748,
\\  N(11)=5662,\, N(13)=5982,\, N(15)=10526,\, N(17)=4028,\, N(19)=11934.
\end{gather*}
\end{remark}

\begin{conjecture} \label{2^k} Any $m\in\Z^+$ can be written as $w^4+x^2+y^2+z^2\ (w,x,y,z\in\N)$
with $x-y\in\{2^k:\ k\in\N\}$.
\end{conjecture}
\begin{remark} See \cite[A350021]{OEIS} for related data. By \cite[Theorem 1.1(ii)]{S19IJNT}, any positive integer can be written as $x^2+y^2+z^2+w^2\ (x,y,z,w\in\N)$ with $x-y\in\{2^a:\ a\in\N\}$.
\end{remark}

For $P(x,y,z,w)\in\Z[x,y,z,w]$, we define its exceptional set
$E(P)$ as the set of all those $n\in\N$ for which there are no $x,y,z,w\in\N$
with $n=x^2+y^2+z^2+w^2$ such that $P(x,y,z,w$ is a square.

\begin{conjecture}\label{134} Any $m\in\N$ not divisible by $8$
can be written as $x^2+y^2+z^2+w^2\ (x,y,z,w\in\N)$ with $x+3y+4z$ a square.
Moreover,
$$E(x+3y+4z)=\{2^{4a+3}q:\ a\in\N,\ q\in\{1,3,5,43\}\}.$$
\end{conjecture}
\begin{remark}\label{R-134} We have verified the former assertion for $m\ls 6\times10^6$.
See \cite[A335624]{OEIS} for related data.
\end{remark}

\begin{conjecture}\label{36} Any $m\in\N$ not divisible by $8$
can be written as $x^2+y^2+z^2+w^2\ (x,y,z,w\in\N)$ with $3x+10y+36z$ a positive square.
Moreover,
$$E(3x+10y+36z)=\{2^{4a+3}q:\ a\in\N,\ q\in\{1,3,5,61\}\}.$$
\end{conjecture}
\begin{remark}\label{R-36} We have verified the former assertion for $m\ls 5\times10^6$. See \cite[A338019]{OEIS} for related data.
\end{remark}

\begin{conjecture}\label{EP} We have
\begin{gather*}
  E(8x+9y)=\{47\times 2^{4a}:\ a\in\N\},\\E(x+4y)=\{2^{4a+2}q:\ a\in\N,\ q\in\{3,23\}\},
\\E(2x+6y+14z)=\{2^{4a+2}q:\ a\in\N,\ q\in\{7,31\}\}.
\end{gather*}
\end{conjecture}
\begin{remark} For other similar conjectures, see \cite[Conjecture 1.20]{S21}.
\end{remark}

\begin{conjecture} Let $a,b,c\in\Z$ with $\gcd(a,b,c)=1$ and $a>0$.
Then any $m\in\N$ can be written as $x^2+y^2+z^2+w^2\ (x,y,z,w\in\Q_{\gs0})$
with $ax+by+cz\in\{4^k:\ k\in\N\}$.
\end{conjecture}
\begin{remark} For example,
$$1261=\l(\f{20}{89}\r)^2+\l(\f{34}{89}\r)^2+\l(\f{35}{89}\r)^2+\l(\f{3160}{89}\r)^2
\ \t{with}\ \f{20}{89}+\f{34}{89}+\f{35}{89}=4^0,$$
$$3=\l(\f1{15}\r)^2+\l(\f 7{15}\r)^2+0^2+\l(\f 53\r)^2
\ \t{with}\ \f1{15}+2\times\f7{15}+4\times0=4^0,$$
and
$$421=\l(\f{38}{105}\r)^2+\l(\f{19}{105}\r)^2+\l(\f{2}{105}\r)^2+\l(\f{718}{35}\r)^2
\ \t{with}\ \f{38}{105}+3\times\f{19}{105}+5\times\f{2}{105}=4^0.$$
\end{remark}

\begin{conjecture} Let $d\in\Z^+$ with $d\not\eq0,1\pmod 4$ and $\gcd(d-1,21)=1$.
Then each $m\in\Z^+$
can be written as $x^2+y^2+z^2+w^2\ (x,y,z,w\in\Q)$ such that $x^2+dy^2\in\{4^k:\ k\in\N\}$.
\end{conjecture}
\begin{remark} For example,
$$643=\l(\f{161}{177}\r)^2+\l(\f{52}{177}\r)^2+\l(\f{49}3\r)^2+\l(\f{1143}{59}\r)^2
$$
with
$$\l(\f{161}{177}\r)^2+2\l(\f{52}{177}\r)^2=4^0=1.$$
For $d=5,15$ we also conjecture that any $m\in\N$ can be written as $x^2+y^2+z^2+w^2\ (x,y,z,w\in\Q)$ such that $x^2+dy^2\in\{4^k:\ k\in\N\}$.
\end{remark}

\begin{conjecture} Let $d\in\Z^+$ with $d\eq3\pmod 8$ and $\gcd(d-1,21)=1$.
If $d\not\in\{3^{4k+1}:\ k\in\N\}$, then each $m\in\Z^+$
can be written as $x^2+y^2+z^2+w^2\ (x,y,z,w\in\Q)$ such that $x^2+dy^2$
is the least power of four not dividing $m$,
\end{conjecture}
\begin{remark} For example,
$$1=\l(\f{10}{27}\r)^2+\l(\f{16}{27}\r)^2+\l(\f 7{27}\r)^2+\l(\f23\r)^2$$
with $$\l(\f{10}{27}\r)^2+11\l(\f{16}{27}\r)^2=4.$$
\end{remark}

\begin{conjecture} If $(a,b,c)$ is among the triples
$$(1,4,7),\ (3,4,5),\ (4,5,11),\ (4,11,13),\ (7,9,12),\ (11,12,13),$$
then any $m\in\N$ can be written as $x^2+y^2+z^2+w^2\ (x,y,z,w\in\Q)$
with $ax^2+by^2+cz^2+dw^2\in\{4^k:\ k\in\N\}$.
\end{conjecture}
\begin{remark} For example,
\begin{gather*}127=9^2+\l(\f{70}{11}\r)^2+\l(\f{15}{11}\r)^2+\l(\f{21}{11}\r)^2
,\ 9^2+4\l(\f{70}{11}\r)^2+7\l(\f{15}{11}\r)^2=4^2;
\\7=\l(\f 53\r)^2+\l(\f{22}{21}\r)^2+\l(\f{17}{21}\r)^2+\l(\f{11}7\r)^2
,\ 3\l(\f 53\r)^2+4\l(\f{22}{21}\r)^2+5\l(\f{17}{21}\r)^2=4^2;
\\7=\l(\f5{18}\r)^2+\l(\f5{18}\r)^2+\l(\f16\r)^2+\l(\f{47}{18}\r)^2,
\ 4\l(\f{5}{18}\r)^2+5\l(\f{5}{18}\r)^2+11\l(\f16\r)^2=4^0;
\\3=\l(\f57\r)^2+\l(\f13\r)^2+\l(\f5{21}\r)^2+\l(\f{32}{21}\r)^2,
\ 4\l(\f 57\r)^2+11\l(\f13\r)^2+13\l(\f 5{21}\r)^2=4;
\\19=\l(\f{23}{49}\r)^2+\l(\f{127}{49}\r)^2+\l(\f{20}{49}\r)^2+\l(\f{169}{49}\r)^2,
\\ 7\l(\f{23}{49}\r)^2+9\l(\f{127}{49}\r)^2+12\l(\f{20}{49}\r)^2=4^3;
\\1=0^2+\l(\f 6{11}\r)^2+\l(\f2{11}\r)^2+\l(\f 9{11}\r)^2,
\ 11\times 0^2+12\l(\f 6{11}\r)^2+13\l(\f2{11}\r)^2=4.
\end{gather*}
\end{remark}

\begin{conjecture} \label{1234} Each $m\in\N$ can be written as $x^2+y^2+z^2+w^2$
with $x,y,z,w\in\N/15$ such that $x^2+2y^2+3z^2+4w^2$ is an integer square.
Moreover, we may even require $x,y,z,w\in\N$ if $m$ is neither divisible by $4$ nor among the following numbers
$$15,\,23,\,26,\,31,\,71,\,77,\,89,\,111,\,127,\,239,\,359,\,575,\,663,\,719,\,991.$$
\end{conjecture}
\begin{remark} For example,
$$26=\l(\f{14}5\r)^2+\l(\f{59}{15}\r)^2+\l(\f{11}{15}\r)^2+\l(\f{22}{15}\r)^2$$
with
$$\l(\f{14}5\r)^2+2\l(\f{59}{15}\r)^2+3\l(\f{11}{15}\r)^2+4\l(\f{22}{15}\r)^2=7^2.$$
\end{remark}

\begin{conjecture} If $(b,c)$ is among the ordered pairs
$$(1,3),\ (1,4),\ (2,2),\ (2,3),\ (3,3),\ (4,4),$$ then any $m\in\N$ can be written as
$x^2+y^2+z^2+w^2$ with $x,y,z,w\in\Q_{\gs0}$ such that $x+by+cz+4w$ is an integer square.
\end{conjecture}
\begin{remark} For example,
\begin{gather*}239=9^2+\l(\f{35}3\r)^2+\l(\f{14}3\r)^2+\l(\f13\r)^2,
\\ 9+\f{35}3+3\times\f{14}3+4\times\f13=6^2;
\\4=\l(\f{15}{13}\r)^2+\l(\f{21}{13}\r)^2+\l(\f1{13}\r)^2+\l(\f 3{13}\r)^2
,\\ \f{15}{13}+\f{21}{13}+4\times\f1{13}+4\times\f3{13}=2^2;
\end{gather*}
\begin{gather*}
8=\l(\f{76}{27}\r)^2+\l(\f 2{27}\r)^2+\l(\f29\r)^2+\l(\f4{27}\r)^2,
\\\f{76}{27}+2\times\f2{27}+2\times\f 29+4\times\f4{27}=2^2;
\\328=\l(\f{124}7\r)^2+2^2+\l(\f47\r)^2+\l(\f{22}7\r)^2,
\\\f{124}7+2\times2+3\times\f47+4\times\f{22}7=6^2;
\\4=\l(\f19\r)^2+\l(\f13\r)^2+\l(\f{17}9\r)^2+\l(\f 59\r)^2,
\\\f19+3\times\f13+3\times\f{17}9+4\times\f 59=3^2;
\\7=\l(\f{41}{17}\r)^2+\l(\f2{17}\r)^2+\l(\f{13}{17}\r)^2+\l(\f{13}{17}\r)^2,
\\\f{41}{17}+4\times\f2{17}+4\times\f{13}{17}+4\times\f{13}{17}=3^2.
\end{gather*}
\end{remark}

\begin{conjecture} [3-5-7 Conjecture] Each $r\in\Q_{\gs0}$ can be written as
$w^2+x^2+y^2+z^2\ (w,x,y,z\in\Q_{\gs0})$ with $3x+5y+7z\in\{a^3:\ a\in\N\}$.
\end{conjecture}
\begin{remark} For example,
$$28432=\l(\f{708}5\r)^2+\l(\f{384}5\r)^2+\l(\f{248}5\r)^2+\l(\f{24}5\r)^2$$
with
$$3\times\f{384}5+5\times\f{248}5+7\times\f{24}5=8^3.$$
We also note that $3x+5y+7z$ in the 3-5-7 conjecture can be replaced by
$x+3y+5z$, but this could be proved by modifying  Machiavelo and N. Tsopanidis's way to show the 1-3-5 conjecture.
\end{remark}

\begin{conjecture} \label{Conj-31} Any $m\in\N$ can be written as $x^2+y^2+z^2+w^2\ (x,y,z,w\in\Q_{\gs0})$
such that $F(x,y)$ is an integer cube, provided that $F(x,y)$ is among the linear polynomials
$$x+6y,\ x+10y,\ x+24y,\ x+31y.$$
\end{conjecture}
\begin{remark} Our computation suggests that $F(x,y)$ also can be replaced by many other linear
polynomials not listed in Conjecture \ref{Conj-31}.
\end{remark}

\begin{conjecture} Each $m\in\N$ can be written as $x^2+y^2+z^2+w^2\ (x,y,z,w\in\Q_{\gs0})$
such that both $x$ and $x+2y+5z$ are integer squares. We may also replace $x+2y+5z$
by any of the following polynomials
\begin{gather*}x+9y,\ x+2y+7z,\ x+y+10z,\ x+3y+10z,\ x+5y+10z,
\\ x+y+4z+8w,\ x+y+7z+9w,\ x+y+7z+10w,
\\ x+2y+7z+10w,\ x+2y+9z+10w,\ x+3y+8z+10w.
\end{gather*}
\end{conjecture}
\begin{remark} Our computation suggests many other such conjectures.
\end{remark}

\begin{conjecture} \label{conj-S} Let $m$ be a positive integer. Then each $n\in\N$ can be written as
$$x^4+y^2+\f{z^4+w^2}m$$ with $x,y,z,w\in\N$, if and only if $m$ belongs to the set
$$\{st^4:\ s\in S \ \t{and}\ t\in\Z^+\},$$ where
\begin{align*}S=&\{5,13,25,65,85,325,18125,20213,32045,68125,
\\&\qquad1105625,1665625,2703125,4250000,5283785\}.
\end{align*}
\end{conjecture}
\begin{remark} For any positive integer $m\ls2\times10^7$, we have verified that
every $n=0,\ldots,10^5$ can be be written as
$x^4+y^2+(z^4+w^2)/m$ with $x,y,z,w\in\N$, if and only if $m$ belongs to the set $\{st^4:\ s\in S
\ \t{and}\ t\in\Z^+\}$
given in Conjecture \ref{conj-S}. See \cite[A349942 and A349945]{OEIS} for related data.
Note that any $r\in\Q_{\gs0}$ can be written as $x^4+y^2+25z^4+w^2$
with $x,y,z,w\in\Q$ if each $n\in\N$ can be written as $x^4+y^2+(z^4+w^2)/25$
with $x,y,z,w\in\Z$.
\end{remark}

\begin{conjecture} \label{8-32} If $(a,b,c,d,m)$ is among the ordered tuples
\begin{align*}&(1,1,7,1,8),\,(1,1,31,1,8),\,(1,1,39,1,8),\,(1,3,7,1,8),\,(2,1,7,1,8),
\\&(2,1,7,1,32),\,(2,1,47,1,8),\,(3,1,7,1,8),\,(3,1,7,1,32),\,(3,1,15,1,8),
\\&(4,1,1,2,3),\,(4,1,1,7,32),\,(4,1,7,1,8),\,(6,1,7,1,8),\,(6,1,23,1,8),
\\&(7,1,2,1,3),\,(13,1,7,1,8)
\end{align*}
or the ordered tuples
\begin{align*}
&(1,1,4,2,3),\,(1,2,1,1,10), \,(1,2,4,1,5),\,(2,1,1,1,5),
\\&(2,1,1,10),\,(2,4,1,1,10),\,(6,1,1,1,5)\,(6,1,1,1,10),\, (6,1,6,1,5),
\end{align*}
then each $n\in\N$ can be written as $ax^4+by^2+(cz^4+dw^2)/m$ with $x,y,z,w\in\N$.
\end{conjecture}
\begin{remark} Conjecture \ref{8-32} were partly announced in \cite[A349945]{OEIS}.
If any $n\in\N$ can be written as
$$2x^4+y^2+\f{7z^4+w^2}{32}=2x^4+y^2+56\l(\f z4\r)^4+2\l(\f w{8}\r)^2$$
with $x,y,z,w\in\N$, then any $r\in\Q_{r\gs0}$ can be written as $2x^4+y^2+56z^4+2w^2$
with $x,y,z,w\in\Q_{\gs0}$. Similar comments work for all other tuples in Conjecture \ref{8-32};
for example, if any $n\in\N$ can be written as
$$x^4+y^2+\f{4z^4+2w^2}{3}=x^4+y^2+108\l(\f z3\r)^4+6\l(\f w{3}\r)^2$$
with $x,y,z,w\in\N$, then any $r\in\Q_{r\gs0}$ can be written as $x^4+y^2+108z^4+6w^2$
with $x,y,z,w\in\Q_{\gs0}$. Note that $108>100$ and this is the reason why we put $(1,1,4,2,3)$
in the second group of ordered tuples in Conjecture \ref{8-32}.
\end{remark}

\begin{conjecture}\label{m^2}
Each $n\in\N$ can be written as $$ax^4 + by^2 + \f{cz^4 + dw^2}{m^2}$$ with $x,y,z,w\in\N$, provided that $(a,b,c,d,m)$ (with $\gcd(c,m) = \gcd(d,m) = 1$ and $a,b,c,d\ls100$) is among the following tuples:
\begin{align*}&(1,1,1,11,50),\,(1,1,1,23,72),\,(1,1,1,31,2000),\,(1,1,1,31,3200),
\\&(1,1,1,47,1176),\,(1,1,2,1,9),\,(1,1,2,7,45),\,(1,1,3,1,2),\,(1,1,3,2,25),
\\&(1,1,7,1,2)\,(1,1,7,1,4),\,(1,1,11,1,2),\,(1,1,15,1,8),\,(1,1,17,7,1800),
\\&(1,1,17,7,2304),\,(1,1,19,1,2),\,(1,1,23,1,8),\,(1,1,23,1,64),\,(1,1,31,1,4),
\\&(1,1,39,1,100),\,(1,1,41,7,324),\,(1,1,47,1,4),\,(1,1,55,1,64),\,(1,1,71,1,16),
\\&(1,1,71,1,36),\,(1,1,79,1,4),\,(1,1,87,1,16),\,(1,1,95,1,4),(1,2,1,1,50),
\\&(1,2,1,7,8),\,(1,2,1,11,50),\,(1,2,1,23,512),\,(1,2,1,47,3528),\,(1,2,3,1,2),
\\&(1,2,3,1,3721),\,(1,2,4,1,25),\,(1,2,5,3,16),\,(1,2,7,1,2),\,(1,2,7,1,4),
\\&(1,2,11,1,9),\,(1,2,11,1,50),\,(1,2,17,7,2304),\,(1,2,19,1,50),\,(1,2,23,1,8),
\\&(1,2,23,1,64),\,(1,2,39,1,64),\,(1,2,47,1,16),\,(1,2,71,1,144),\,(1,2,79,1,2704),
\\&(1,3,1,1,5),\,(1,3,2,1,729),\,(1,3,7,1,4),\,(1,3,7,1,32),\,(1,3,11,1,450),
\\&(1,3,23,1,32),\,(1,3,23,1,576),\,(1,3,31,1,1600),\,(1,3,47,1,16),\,(1,3,47,1,36),
\\&(1,3,71,1,400),\,(1,4,2,1,9),\,(1,4,71,1,144),\,(1,4,79,1,400),\,(1,6,31,1,1600),
\\&(1,7,71,1,3600),\,(1,8,47,1,2304),\,(1,11,23,1,576),
\end{align*}
\begin{align*}&(2,1,1,1,5),\,(2,1,1,1,10),\,(2,1,1,1,25),\,(2,1,1,4,125),\,(2,1,1,11,150),
\\&(2,1,1,11,450),\,(2,1,1,19,1250),\,(2,1,1,23,24),\,(2,1,1,23,144),
\\&(2,1,1,23,1152),\,(2,1,1,31,800),\,(2,1,1,31,1280),\,(2,1,2,1,3),\,(2,1,3,1,2),
\\&(2,1,4,1,5),\,(2,1,5,1,3),\,(2,1,5,19,432),\,(2,1,11,1,2),\,(2,1,11,1,225),
\\&(2,1,15,1,8),\,(2,1,17,7,144),\,(2,1,17,7,200),\,(2,1,19,1,50),\,(2,1,23,1,8),
\\&(2,1,23,1,16),\,(2,1,29,3,3136),\,(2,1,31,1,4),(2,1,39,1,100),\,(2,1,41,7,324),
\\&(2,1,47,1,4),\,(2,1,71,1,4),\,(2,1,79,1,16),\,(2,1,95,1,4),\,(2,2,1,7,1936),
\\&(2,2,7,1,16),\,(2,2,23,1,4608),\,(2,2,39,1,400),\,(2,2,47,1,144),(2,2,71,1,144),
\\&(2,2,79,1,400),\,(2,3,1,1,4225),\,(2,3,7,1,8),\,(2,3,23,1,32),\,(2,3,23,1,36),
\\&(2,3,47,1,1764),\,(2,3,71,1,1600),\,(2,4,1,1,25),\,(2,4,1,1,845),\,(2,4,23,1,36),
\\&(2,4,23,1,288),\,(2,4,47,1,16),\,(2,4,71,1,900),
\end{align*}
\begin{align*}&(3,1,1,2,3),\,(3,1,1,2,9),\,(3,1,1,23,648),\,(3,1,2,7,45),\,(3,1,7,1,8),
\\&(3,1,7,1,16),\,(3,1,11,1,50),\,(3,1,19,1,1250),\,(3,1,23,1,18),\,(3,1,23,1,36),
\\&\,(3,1,47,1,1764),\,(3,1,71,1,144),\,(3,1,79,1,400),\,(3,1,95,1,144),\,(3,2,1,1,125),
\\&(3,2,1,1,169),\,(3,2,1,23,576),(3,2,3,1,98),\,(3,2,7,1,4),\,(3,2,7,1,8),
\\&(3,2,11,1,450),\,(3,2,23,1,16),\,(3,2,23,1,32),\,(3,2,71,1,144),
\\&(3,2,79,1,100),\,(3,3,23,1,288),\,(3,3,71,1,3600),\,(3,4,71,1,3600),
\\&(4,1,1,2,27),\,(4,1,2,1,3),\,(4,1,2,1,81),\,(4,1,3,2,25),\,(4,1,7,1,8),
\\&(4,1,7,1,16),\,(4,1,11,1,450),\,(4,1,23,1,8),\,(4,1,23,1,36),\,(4,1,31,1,16),
\\&(4,1,47,1,36),\,(4,1,71,1,144),\,(4,1,79,1,400),\,(4,2,1,1,5),\,(4,2,1,23,576),
\\&(4,2,7,1,2),\,(4,2,7,1,4),\,(4,2,11,1,25),\,(4,2,23,1,16),\,(4,2,23,1,32),
\\&(4,2,47,1,16),\,(4,2,47,1,36),\,(4,2,71,1,144)\,(4,2,79,1,1600),\,(4,2,87,1,16),
\\&(4,3,7,1,32),\,(4,3,23,1,576),\,(4,3,71,1,3600),\,(4,4,71,1,3600),
\end{align*}
\begin{align*}&(5,1,1,2,27),\,(5,1,2,1,3),\,(5,1,2,1,9),\,(5,1,7,1,2),\,(5,1,23,1,16),
\\&(5,1,23,1,18),\,(5,1,31,1,16),\,(5,1,39,1,4),\,(5,1,47,1,36),\,(5,1,71,1,16),
\\&(5,2,1,1,3125),\,(5,2,1,7,16),\,(5,2,1,23,144),\,(5,2,11,1,450),
\\&(5,2,23,1,2592),\,(5,2,31,1,196),\,(5,2,47,1,2304),\,(5,3,7,1,3872),
\\&(5,3,71,1,3600),\,(6,1,1,7,32),\,(6,1,1,11,4050),\,(6,1,1,23,192),
\\&(6,1,1,23,1152),\,(6,1,7,1,8),\,(6,1,7,1,484)\,(6,1,23,1,8),\,(6,1,23,1,36),
\\&(6,1,31,1,400),\,(6,1,47,1,36),\,(6,1,71,1,144),\,(6,2,47,1,2304),
\\&(6,2,71,1,3600),\,(6,4,71,1,3600),\,(7,1,11,1,50),\,(7,1,23,1,18),
\\&(7,1,23,1,576),\,(7,1,47,1,144),\,(7,1,71,1,400),\,(7,2,11,1,18),
\\&(7,2,11,1,50),\,(7,2,23,1,32),\,(7,2,23,1,72),\,(7,2,23,1,576),
\\&(7,2,31,1,100),\,(7,2,47,1,2304),\,(7,3,47,1,2304),\,(7,3,71,1,3600),
\end{align*}
\begin{align*}
&(8,1,1,23,288),\,(8,1,7,1,1936),\,(8,1,11,1,450),\,(8,1,23,1,32),
\\&(8,1,23,1,576),\,(8,1,31,1,16),\,(8,1,47,1,2304),\,(8,1,71,1,144),
\\&(8,2,47,1,1296),\,(9,1,7,1,8),\,(9,1,11,1,450),\,(9,1,23,1,36),
\\&(9,1,23,1,288),\,(9,1,47,1,1764),\,(9,1,71,1,400),\,(9,2,3,1,2),
\\&(9,2,31,1,400),\,(9,2,47,1,1764),\,(9,2,71,1,144),\,(9,3,23,1,1296),
\\&(10,1,7,1,8),\,(10,1,7,1,16),\,(10,1,11,1,450),\,(10,1,23,1,4),
\\&(10,1,23,1,8),\,(10,1,47,1,36),\,(10,1,71,1,400),\,(10,2,31,1,400),
\\&(10,2,47,1,576),\,(10,3,71,1,3600),
\end{align*}
\begin{align*}
&(11,1,1,2,81),\,(11,1,2,1,9),\,(11,1,2,1,27),\,(11,1,4,2,27),
\\&(11,1,7,1,8),\,(11,1,23,1,18),\,(11,1,23,1,36),\,(11,1,31,1,100),
\\&(11,1,47,1,324),\,(11,1,71,1,3600),\,(11,2,31,1,16),\,(11,2,47,1,144),
\\&(11,2,71,1,144),\,(11,2,79,1,100),\,(11,3,47,1,2304),\,(11,3,71,1,3600),
\\&(12,1,2,1,27),\,(12,1,2,1,81),\,(12,1,7,1,8),\,(12,1,23,1,32),\,(12,1,23,1,36),
\\&(12,1,31,1,400),\,(12,1,47,1,144),\,(12,1,71,1,400),\,(12,2,47,1,576),
\\&(12,2,71,1,144),\,(13,1,1,2,9),\,(13,1,2,1,27),\,(13,1,7,1,16),
\\&(13,1,11,1,450),\,(13,1,23,1,32),\,(13,1,47,1,2304),\,(13,1,71,1,3600),
\\&(13,2,31,1,1600),\,(13,2,47,1,576),\,(13,2,71,1,144),\,(13,3,47,1,2304),
\\&(14,1,11,1,50),\,(14,1,23,1,4),\,(14,1,23,1,18),\,(14,1,47,1,144),
\\&(14,1,71,1,900),\,(15,1,23,1,18),\,(16,1,23,1,72),\,(16,1,71,1,576),
\\&(17,1,7,1,8),\,(17,2,71,1,144),\,(18,1,7,1,8),\,(18,1,11,1,450),
\\&(18,1,23,1,72),\,(18,1,71,1,576),\,(19,1,7,1,968),\,(19,1,11,1,450),
\\&(19,1,71,1,64),\,(19,2,71,1,3600),,\,(20,1,11,1,450),\,(20,1,23,1,72),
\\&(20,1,71,1,400),\,(21,1,23,1,72),\,(21,1,71,1,3600),\,(21,2,71,1,3600),
\\&(22,1,23,1,1152),\,(22,1,71,1,400),\,(22,2,71,3600),(23,1,71,1,400),
\\&(23,2,71,1,900),\,(24,1,23,1,1152),\,(24,1,71,1,3600),\,(25,1,23,1,288),
\\&(25,1,71,1,900),\,(26,1,23,1,288),\,(26,2,71,1,3600),\,(28,1,23,1,1152),
\\&(28,1,71,1,900),\,(29,1,71,1,900),\,(30,1,23,1,72),\,(31,1,71,1,900),
\\&(32,1,23,1,72),\,(33,1,71,1,3600),\,(34,1,23,1,72),\,(34,1,71,1,900),
\\&(35,1,71,1,900),\,(36,1,23,1,288),\,(36,1,71,1,3600),\,(37,1,71,1,900),
\\&(38,1,71,1,900),\,(39,1,23,1,288),\,(40,1,23,1,288),\,(40,1,71,1,3600),
\\&(43,1,23,1,288),\,(43,1,71,1,3600),\,(44,1,23,1,72),\,(44,1,71,1,3600),
\\&(45,1,23,1,288),\,(46,1,71,1,3600),\,(47,1,23,1,288),\,(47,1,71,1,3600),
\\&(48,1,23,1,288),\,(49,1,23,1,288),\,(49,1,71,1,3600),
\end{align*}
\begin{align*}
&(51,1,71,1,3600),\,(52,1,23,1,288),\,(52,1,71,1,3600),\,(53,1,23,1,72),
\\&(53,1,71,1,3600),\,(54,1,71,1,3600)\,(55,1,23,1,288),\,(55,1,71,1,3600),
\\&(56,1,23,1,72),\,(57,1,23,1,72),\,(58,1,71,1,3600),\,(59,1,71,1,3600),
\\&(60,1,71,1,3600),\,(62,1,23,1,288),\,(63,1,23,1,288),\,(66,1,23,1,288),
\\&(67,1,23,1,288),\,(68,1,23,1,288).
\end{align*}
\end{conjecture}
\begin{remark} We obtain the tuples in Conjecture \ref{m^2}
via a computer search. For each tuple $(a,b,c,d,m)$ listed in Conjecture \ref{m^2}, we have
verified that every $n=0,\ldots,30000$ can be written as $ax^4+by^2+(cz^4+dw^2)/m^2$
with $x,y,z,w\in\N$.
\end{remark}

Conjectures \ref{conj-S}--\ref{m^2} implies that any $n\in\N$ can be written
as $ax^4+by^2+cz^4+dw^2$ with $x,y\in\N$ and $z,w\in\Q_{\gs0}$, provided that
$(a,b,c,d,m)$ with $a,b,c,d\ls100$ is among the following tuples:
\begin{align*}&(1,1,2,1),\,(1,1,3,2),\,(1,1,4,11),\,(1,1,4,23),\,(1,1,4,31),\,(1,1,7,1),
\\&(1,1,12,1),\,(1,1,14,2),\,(1,1,17,7),\,(1,1,23,1),\,(1,1,25,1),\,(1,1,25,31),
\\&(1,1,28,1),(1,1,31,1),\,(1,1,36,47),\,(1,1,39,1),\,(1,1,41,7),\,(1,1,44,1),
\\&(1,1,47,1),\,(1,1,50,7),\,(1,1,55,1),\,(1,1,60,1),\,(1,1,62,2),\,(1,1,68,7),
\\&(1,1,71,1),\,(1,1,76,1),\,(1,1,78,2),\,(1,1,79,1),\,(1,1,87,1),\,(1,1,92,1),
\\&(1,1,95,1),\,(1,2,3,1),\,(1,2,4,1),\,(1,2,4,7),\,(1,2,4,11),\,(1,2,4,23),
\\&(1,2,4,47),\,(1,2,5,3),\,(1,2,7,1),\,(1,2,11,1),\,(1,2,12,1),\,(1,2,17,7),
\\&(1,2,23,1),\,(1,2,28,1),\,(1,2,39,1),\,(1,2,44,1),\,(1,2,47,1),\,(1,2,76,1),
\\&(1,2,71,1),\,(1,2,79,1),\,(1,2,92,1),\,(1,3,2,1),\,(1,3,7,1),\,(1,3,14,2),
\\&(1,3,23,1),\,(1,3,25,1),\,(1,3,28,1),\,(1,3,31,1),\,(1,3,44,1),\,(1,3,47,1),
\\&(1,3,71,1),\,(1,3,92,1),\,(1,4,2,1),\,(1,4,71,1),\,(1,4,79,1),\,(1,6,31,1),
\\&(1,7,71,1),\,(1,8,47,1),\,(1,11,23,1),
\end{align*}
\begin{align*}&(2,1,1,1),\,(2,1,1,23),\,(2,1,4,11),\,(2,1,4,19),\,(2,1,4,23),\,(2,1,4,31),
\\&(2,1,11,1),\,(2,1,12,1),\,(2,1,14,2),\,(2,1,17,7),\,(2,1,18,1),\,(2,1,23,1),
\\&(2,1,25,1),\,(2,1,25,31),\,(2,1,29,3),\,(2,1,31,1),\,(2,1,36,11),\,(2,1,36,23),
\\&(2,1,39,1),\,(2,1,41,7),\,(2,1,44,1),\,(2,1,45,1),\,(2,1,45,19),\,(2,1,47,1),
\\&(2,1,56,2),\,(2,1,60,1),\,(2,1,68,7),\,(2,1,71,1)\,(2,1,76,1),\,(2,1,79,1),
\\&(2,1,92,1),\,(2,1,94,2),\,(2,1,95,1),\,(2,1,100,1),\,(2,2,1,7),\,(2,2,7,1),
\\&(2,2,39,1),\,(2,2,47,1),\,(2,2,71,1),\,(2,2,79,1),\,(2,2,92,1),\,(2,3,1,1),
\\&(2,3,23,1),\,(2,3,28,1),\,(2,3,47,1),\,(2,3,71,1),\,(2,3,92,1),\,(2,4,1,1),
\\&(2,4,23,1),\,(2,4,25,1),\,\,(2,4,47,1),\,(2,4,71,1),\,(2,4,92,1),
\end{align*}
\begin{align*}&(3,1,1,2),\,(3,1,1,92),\,(3,1,7,1),\,(3,1,9,2),\,(3,1,14,2),\,(3,1,23,1),
\\&(3,1,28,1),\,(3,1,30,2),\,(3,1,44,1),\,(3,1,47,1),\,(3,1,50,7),\,(3,1,56,2),
\\&(3,1,71,1),\,(3,1,76,1),\,(3,1,79,1),\,(3,1,92,1),\,(3,1,95,1),\,(3,2,1,1),
\\&(3,2,1,23),\,(3,2,7,1),\,(3,2,12,1),\,(3,2,23,1),\,(3,2,25,1),\,(3,2,28,1),
\\&(3,2,44,1),\,(3,2,71,1),\,(3,2,79,1),\,(3,2,92,1),\,(3,3,71,1),\,(3,3,92,1),
\\&(3,4,71,1),\,(4,1,2,1),\,(4,1,3,2),\,(4,1,7,1),\,(4,1,8,14),\,(4,1,9,2),
\\&(4,1,14,2),\,(4,1,18,1),\,(4,1,23,1),\,(4,1,27,6),\,(4,1,28,1),\,(4,1,31,1),
\\&(4,1,44,1),\,(4,1,47,1),\,(4,1,71,1),\,(4,1,79,1),\,(4,1,92,1),\,(4,2,1,1,5),
\\&(4,2,1,23),\,(4,2,7,1),\,(4,2,11,1),\,(4,2,23,1)\,(4,2,28,1),\,(4,2,47,1),
\\&(4,2,71,1)\,(4,2,79,1),\,(4,2,87,1),\,(4,2,92,1),\,(4,3,23,1),\,(4,3,28,1),
\\&(4,3,71,1),\,(4,4,71,1),\,(5,1,2,1),\,(5,1,9,2),\,(5,1,18,1),\,(5,1,23,1),
\\&(5,1,28,1),\,(5,1,31,1),\,(5,1,39,1),\,(5,1,47,1),\,(5,1,71,1),\,(5,1,92,1),
\\&(5,2,1,7),\,(5,2,1,23),\,(5,2,25,1),\,(5,2,31,1),\,(5,2,44,1),\,(5,2,47,1),
\\&(5,2,92,1),\,(5,3,28,1),\,(5,3,71,1),\,(6,1,4,7),\,(6,1,4,11),\,(6,1,4,23),
\\&(6,1,7,1),\,(6,1,9,23),\,(6,1,14,2),\,(6,1,23,1),\,(6,1,28,1),\,(6,1,31,1),
\\&(6,1,46,2),\,(6,1,47,1),\,(6,1,71,1),\,(6,1,92,1),\,(6,2,47,1),\,(6,2,71,1),
\\&(6,4,71,1),\,(7,1,23,1),\,(7,1,44,1),\,(7,1,47,1),\,(7,1,54,3),\,(7,1,71,1),
\\&(7,1,92,1),\,(7,2,23,1),\,(7,2,31,1),\,(7,2,44,1),(7,2,47,1),\,(7,2,92,1),
\\&(7,3,47,1),\,(7,3,71,1),\,(8,1,4,23),\,(8,1,7,1),(8,1,23,1),\,(8,1,31,1),
\\&(8,1,44,1),\,(8,1,47,1),\,(8,1,71,1),\,(8,1,92,1)\,(8,2,47,1),
\\&(9,1,23,1),\,(9,1,28,1),\,(9,1,44,1),\,(9,1,47,1),\,(9,1,71,1),\,(9,1,92,1),
\\&(9,2,12,1),\,(9,2,31,1),\,(9,2,47,1),\,(9,2,71,1),\,(9,3,23,1),\,(10,1,7,1),
\\&(10,1,23,1),\,(10,1,28,1),\,(10,1,44,1),\,(10,1,47,1),\,(10,1,71,1),(10,1,92,1),
\\&(10,2,31,1),\,(10,2,47,1),\,(10,3,71,1),
\end{align*}
\begin{align*}
&(11,1,1,2),\,(11,1,2,1),\,(11,1,18,1),\,(11,1,23,1),\,(11,1,28,1),\,(11,1,31,1),
\\&(11,1,36,2),\,(11,1,47,1),\,(11,1,71,1),\,(11,1,92,1),\,(11,2,31,1),\,(11,2,47,1),
\\&(11,2,71,1),\,(11,2,79,1),\,(11,3,47,1),\,(11,3,71,1),\,(12,1,2,1),\,(12,1,18,1),
\\&(12,1,23,1),\,(12,1,28,1),(12,1,31,1),\,(12,1,47,1),\,(12,1,71,1),\,(12,1,92,1),
\\&(12,2,47,1),\,(12,2,71,1),\,(13,1,1,2),\,(13,1,7,1),\,(13,1,14,2),(13,1,18,1),
\\&(13,1,44,1),\,(13,1,47,1),\,(13,1,71,1),\,(13,1,92,1),\,(13,2,31,1),\,(13,2,47,1),
\\&(13,2,71,1),\,(13,3,47,1),\,(14,1,23,1),\,(14,1,44,1),\,(14,1,47,1),\,(14,1,71,1),
\\&(14,1,92,1),\,(15,1,92,1),\,(16,1,71,1),\,(16,1,92,1),\,(17,1,28,1),\,(17,2,71,1),
\\&(18,1,28,1),\,(18,1,44,1),\,(18,1,71,1),\,(18,1,92,1),\,(19,1,28,1),\,(19,1,44,1),
\\&(19,1,71,1),\,(19,2,71,1),\,(20,1,14,2),\,(20,1,44,1),\,(20,1,71,1),\,(20,1,92,1),
\end{align*}

\begin{align*}
&(21,1,71,1),\,(21,1,92,1),\,(21,2,71,1),\,(22,1,71,1),\,(22,1,92,1),\,(22,2,71,1),
\\&(23,1,71,1),\,(23,2,71,1),\,(24,1,71,1),\,(24,1,92,1),\,(25,1,71,1),\,(25,1,92,1),
\\&(26,1,92,1),\,(26,2,71,1),\,(28,1,71,1),\,(28,1,92,1),\,(29,1,71,1),\,(30,1,92,1),
\\&(31,1,71,1),\,(32,1,92,1),\,(33,1,71,1),\,(34,1,71,1),\,(34,1,92,1),\,(35,1,71,1),
\\&(36,1,71,1),\,(36,1,92,1),\,(37,1,71,1),\,(38,1,71,1),\,(39,1,92,1),\,(40,1,71,1),
\\&(40,1,92,1),\,(43,1,71,1),\,(43,1,92,1),\,(44,1,71,1),\,(44,1,92,1),\,(45,1,92,1),
\\&(46,1,71,1),\,(47,1,71,1),\,(47,1,92,1),\,(48,1,92,1),\,(49,1,71,1),\,(49,1,92,1),
\\&(51,1,71,1),\,(52,1,71,1),\,(52,1,92,1),\,(53,1,71,1),\,(53,1,92,1),\,(54,1,71,1)
\\&(55,1,71,1),\,(55,1,92,1),\,(56,1,92,1),\,(57,1,92,1),\,(58,1,71,1),\,(59,1,71,1),
\\&(60,1,71,1),\,(62,1,92,1),\,(63,1,92,1),\,(66,1,92,1),(67,1,92,1),\,(68,1,92,1).
\end{align*}

\begin{conjecture}[23-Conjecture] Each $n\in\N$ can be written as $w^4+x^4+(y^2+23z^2)/16$,
where $w$ is zero or a power of two (including $2^0=1$), and $x,y,z$ are nonnegative integers.
\end{conjecture}
\begin{remark} We have verified this for $n\ls 3\times10^6$. See \cite[A347827]{OEIS} for related data.
We also conjecture that any positive integer can be written as $w^4+x^2+(y^2+23z^4)/324$
with $w\in\{2^k:\ k\in\N\}$ and $x,y,z\in\N$, see \cite[A347562]{OEIS} for related data.
\end{remark}

\begin{conjecture} \label{4-9} {\rm (i)} Any $n\in\N$
can be written as $w^4+(x^4+y^2+z^2)/9$ with $w,x,y,z\in\N$. 
Also, each $n\in\N$ can be written as $w^4+(4x^4+y^2+z^2)/81$ with $w,x,y,z\in\N$.

{\rm (ii)} Each $n\in\N$ can be written as $w^8+(x^4+y^2+z^2)/18$ with $w,x,y,z\in\N$.
Also, any $n\in\N$ can be written as $w^6+(x^6+y^2+z^2)/882$ with $w,x,y,z\in\N$.
\end{conjecture}
\begin{remark} See \cite[A349943]{OEIS} for related data.
Conjecture \ref{4-9}(i) implies that for $d=2,3$ each $r\in\Q_{\gs0}$
can be written as $w^4+d^2x^4+y^2+z^2$ with $w,x,y,z\in\Q$.
We also conjecture that for $d=3,5$ any $n\in\N$ can be written as $2w^4+(x^4+y^2+z^2)/d^2$ with $w,x,y,z\in\N$.
\end{remark}

\begin{conjecture} \label{a^2} {\rm (i)} For any $a\in\Z^+$, there is a positive integer $m$ such that each $n\in\N$
can be written as $w^4+(a^2x^4+y^2+z^2)/m^4$ with $w,x,y,z\in\N$. 

{\rm (ii)} For any positive odd number $a$, there is a positive integer $m$ such that each $n\in\N$
can be written as $4w^4+(a^2x^4+y^2+z^2)/m^4$ with $w,x,y,z\in\N$. 
\end{conjecture}
\begin{remark} We have checked this for $a\ls 100$. For example, in the case $a=1$ we guess that $m=3$
works for both parts of Conjecture \ref{a^2} (cf. \cite[A350860]{OEIS}). 
Conjecture \ref{a^2}(i) implies that for any $a\in\Z^+$ each $r\in \Q_{\gs0}$ can be written as
$w^4+a^2x^4+y^2+z^2$ with $w,x,y,z\in\Q_{\gs0}$.
\end{remark}

\begin{conjecture} {\rm (i)} Any $n\in\Z^+$ can be written as $x^4+y^2+(z^2+2^{4\da}11)/60$
with $x,y,z\in\N$ and $\da\in\{0,1\}$. Also, each $n\in\Z^+$ can be written as $2x^4+y^2+(z^2+11\times2^{4\da})/60$ with $x,y,z\in\N$ and $\da\in\{0,1\}$.

{\rm (ii)} Each $n\in\Z^+$ can be written as $ax^4+y^2+(z^2+b2^{4w})/m$
with $x,y,z,w\in\N$ provided that $(a,b,m)$ is among the following ordered triples:
$$(1,11,15),\ (1,15,64),\ (1,23,24),\ (2,11,15),\ (2,14,15),\ (2,23,24),\ (5,11,60).$$
\end{conjecture}
\begin{remark} See \cite[A349957]{OEIS} for related data.
\end{remark}

\begin{conjecture} {\rm (i)} Any $n\in\Z^+$ can be written as $x^4+y^2+(z^2+2^{2\da+1})/3$
with $x,y,z\in\N$ and $\da\in\{0,1\}$.

{\rm (ii)} Any $n\in\Z^+$ can be written as $x^4+y^2+(z^2+a4^{\da})/m$
with $x,y,z\in\N$ and $\da\in\{0,1\}$, provided that $(a,m)$ is among the following ordered pairs:
$$(11,12),\ (11,60),\ (14,15),\ (23,24),\ (23,32),\ (23,48).$$

{\rm (iii)} Any $n\in\Z^+$ can be written as $F(x,y,z,\da)$ with $x,y,z\in\N$ and $\da\in\{0,1\}$
provided that $F(x,y,z,\da)$ is among the following expressions:
\begin{gather*}x^4+2y^2+\f{z^2+23\times4^{\da}}{96},\ 2x^4+y^2+\f{z^2+11\times4^{\da}}{60},
\\ 2x^4+y^2+\f{z^2+23\times4^{\da}}{24},
\ 2x^4+y^2+\f{z^2+23\times4^{\da}}{48},
\ 4x^4+y^2+\f{z^2+23\times4^{\da}}{48}.
\end{gather*}
\end{conjecture}
\begin{remark} See \cite[A349992]{OEIS} for related data. For example,
both $80$ and $38064$ have a unique representation required in part (i):
$$80=1^4+6^2+\f{11^2+2^3}3\ \ \t{and}\ \ 38064=3^4+157^2+\f{200^2+2}3.$$
\end{remark}

Conjecture \ref{2^k} implies that each $m\in\Z^+$ can be written as $x^4+y^2+(z^2+4^w)/2$ with $x,y,z,w\in\N$ (cf. \cite[A349661]{OEIS}). In contrast, we have the following conjecture.

\begin{conjecture} {\rm (i)} Any $n\in\Z^+$ can be written as $4x^4+y^2+(z^2+4^w)/2$ with $x,y,z,w\in\N$.

{\rm (ii)} Each $n\in\Z^+$ can be written as $ax^4+by^2+(cz^2+d4^w)/m$
with $x,y,z,w\in\N$, provided that $(a,b,c,d,m)$ is among the following ordered tuples:
\begin{align*}&(1,1,1,2,51),\, (1,1,2,3,5),\, (1,1,1,6,7),\, (1,1,2,7,15),\,(1,1,1,11,20),
\\&(1,1,1,11,36),\,(1,1,2,13,15),\,(1,1,1,15,16),\,(1,1,1,17,66),\,(1,1,1,19,20),
\\&(1,2,1,1,5),\,(1,2,1,2,3),\,(1,2,1,3,4),\,(1,2,1,7,8),\,(1,2,1,7,16),
\\&(1,2,1,11,12),\,(1,2,1,11,36),\,(1,2,1,11,60),\,(1,2,1,17,21),\,(2,1,1,1,5),
\\&(2,1,2,1,3),\,(2,1,1,3,4),\,(2,1,2,3,5),\,(2,1,1,5,9),\,(2,1,1,5,21),
\\&(2,1,1,6,7),\,(2,1,1,11,12),\,\,(2,1,2,7,15),\,(3,1,1,2,3),\,(3,1,1,7,8),
\\&(4,1,2,1,3),\,(4,1,1,2,8),\,(4,2,1,2,3),\,(6,1,1,1,5),\,(9,1,1,2,3).
\end{align*}
\end{conjecture}
\begin{remark} Part (i) is a new refinement of Lagrange's four-square theorem since
$$\f{x^2+y^2}2=\l(\f{x+y}2\r)^2+\l(\f{x-y}2\r)^2.$$
See \cite[A350012]{OEIS} for related data.
\end{remark}

\begin{conjecture} Let $n\in\Z^+$ with $8\nmid n$ and $n\not\in\{4^k(8l+7):\ k,l\in\N\}$.
Then we can write $n=x^2+y^2+z^2$ with $x,y,z\in\Q_{\gs0}$ such that $x+2y+3z\in\Q^2$.
Also, we can write $n=x^2+y^2+z^2$ with $x,y,z\in\Q_{\gs0}$ such that $x+y+6z\in\Q^2$.
\end{conjecture}
\begin{remark} The author \cite[Conjecture 4.2(ii)]{S19IJNT} conjectured that
any $n\in\N\sm\{4^k(8l+7):\ k,l\in\N\}$ can be written as $x^2+y^2+z^2$ with $x,y,z\in\Z$
such that $x+2y+3z$ is a square or twice a square.
\end{remark}

It is known (cf. \cite[pp.\,112-113]{D39}) that for each $d=1,3,4$ any positive odd integer can be written as
$x^2+2y^2+dz^2$ with $x,y,z\in\N$.

\begin{conjecture} Let $d\in\{1,3,4\}$. Then any positive odd integer can be written as
$x^2+2y^2+dz^2$ with $x,y,z\in\Q_{\gs0}$ such that $x+y+z\in\Q^2$.
\end{conjecture}
\begin{remark}  For example,
\begin{align*}3=&\l(\f{23}{25}\r)^2+\l(\f{36}{25}\r)^2+2\l(\f15\r)^2
\ \ \t{and}\ \ \f{23}{25}+\f{36}{25}+\f15=\l(\f 85\r)^2,
\\7=&\l(\f{22}{27}\r)^2+2\l(\f{16}{27}\r)^2+3\l(\f{37}{27}\r)^2\ \ \t{and}\
\ \f{22}{27}+\f{16}{27}+\f{37}{27}=\l(\f 53\r)^2,
\\9=&\l(\f3{10}\r)^2+2\l(\f12\r)^2+4\l(\f{29}{20}\r)^2\
\ \t{and}\ \ \f3{10}+\f12+\f{29}{20}=\l(\f 32\r)^2.
\end{align*}
\end{remark}

\begin{conjecture} Each $r\in\Q_{\gs0}$ can be written as $x^7+y^2+z^2$
with $x,y,z\in\Q_{\gs0}$. Moreover, for any $n\in\N$ we can write $6^{14}n$
as $x^7+y^2+z^2$ with $x,y,z\in\N$.
\end{conjecture}
\begin{remark} For any $n\in\N$, we also conjecture that
$2^6n$ can be written as $x^3+y^2+z^2$ with $x,y,z\in\N$, and that $2^{20}n$
can be written as $x^5+y^2+z^2$ with $x,y,z\in\N$.
\end{remark}

\begin{conjecture} For any $n\in\N$, we can write $6552n$ as $w^6+x^3+y^3+z^3$
with $w,x,y,z\in\N$.
\end{conjecture}
\begin{remark} We also conjecture that for any $n\in\N$ we can write $504n$ as $w^3+x^3+y^3+z^3$
with $w,x,y,z\in\N$.
\end{remark}


\end{document}